\documentclass[12pt]{amsart}
\usepackage{amssymb}
\usepackage{amsfonts}
\usepackage{amscd}
\usepackage{latexsym}
\usepackage{enumerate}
\usepackage{amsmath}
\usepackage[all]{xy}
\numberwithin{equation}{section}

\newcounter{cs}
\stepcounter{cs}
\newcounter{ds}
\stepcounter{ds}

\newcommand{\casos}{\begin{itemize}}
\newcommand{\fcasos}{\end{itemize}\setcounter{cs}{1}}

\newcommand{\ol}{\overline}

\newcommand{\im}{\mbox{Im}}

\newcommand{\Bla}{{\mathcal Bla}}

\newcommand{\Si}{{\rm Sink}}

\newtheorem{lem}{Lemma}[section]
\newtheorem{cor}[lem]{Corollary}
\newtheorem{theor}[lem]{Theorem}
\newtheorem{prop}[lem]{Proposition}

\theoremstyle{definition}
\newtheorem{defi}[lem]{Definition}

\newtheorem{exems}[lem]{Examples}

\newtheorem{remark}[lem]{Remark}

\newcommand{\N}{\mathbb{N}}

\newcommand{\Z}{\mathbb{Z}}

     \newcommand{\id}[1][]{\mathbf{1}_{#1}}             %Identitat
     \newcommand{\supp}[1]{\mathrm{supp}(#1)}

     \DeclareMathOperator{\coker}{coker}

     \newcommand{\cb}[0]{k}                             %Cos base

     \newcommand{\fp}[1]{\text{\bf fp}(#1)}
     \newcommand{\fpl}[1]{\text{\bf fp}(#1)_{\mathrm{fl}}}
     \newcommand{\fnp}[1]{\text{\bf fnp}(#1)_{\mathrm{fl}}}
     \newcommand{\Mod}[1]{#1\text{\bf -Mod}}

     \newcommand{\cam}[2][]{P_{#1}(#2)}
     \newcommand{\rat}[1]{P_{{\rm rat}}(#1)}
     \newcommand{\ser}[1]{P((#1))}
     \newcommand{\reg}[1]{Q(#1)}
     \newcommand{\lev}[2][]{L_{#1}(#2)}

\begin{document}
\title[Module theory over Leavitt path algebras]{Module theory over Leavitt
path algebras\\ and $K$-theory}
\author{Pere Ara and Miquel Brustenga}\address{Departament de Matem\`atiques, Universitat Aut\`onoma de
Barcelona, 08193, Bellaterra (Barcelona), Spain}\email{para@mat.uab.cat, mbrusten@mat.uab.cat}

\thanks{Both
authors were partially supported by DGI MICIIN-FEDER
MTM2008-06201-C02-01, and by the Comissionat per Universitats i
Recerca de la Generalitat de Catalunya. The second author was
partially
 supported by a grant of the Departament de Matem\`{a}tiques,
 Universitat Aut\`{o}noma de Barcelona.}
 \date{\today}
 \begin{abstract} Let $k$ be a field and let $E$ be a finite quiver.
 We study the structure of the finitely presented modules of finite length over the
 Leavitt path algebra $L_k (E) $ and show its close relationship
 with the finite-dimensional representations of the inverse quiver $\overline{E}$
 of $E$, as well as with the class of finitely generated
 $P _k(E)$-modules $M$ such that ${\rm Tor}_q^{P_k (E)}(k^{|E^0|},M)=0$ for all $q$,
 where $P_k(E)$ is the usual path algebra of $E$.
By using these results we compute the higher $K$-theory of the von
Neumann regular algebra $Q_k (E)=L_k (E)\Sigma^{-1}$, where $\Sigma
$ is the set of all square matrices over $P_k (E)$ which are sent to
invertible matrices by the augmentation map $\epsilon \colon P_k
(E)\to k^{|E^0|}$.
 \end{abstract}
\maketitle

\section{Introduction}

\noindent For a field $\cb $ and an integer $n\ge 2$, the Leavitt
algebra $L(1,n)$ of type $(1,n)$ is the algebra with generators
$x_i,y_j$, $1\le i,j\le n$ with defining relations given by
$$(x_1,\dots ,x_n)(y_1,\dots ,y_n)^t =1,\qquad (y_1,\dots ,y_n)^t(x_1,\dots ,x_n)=I_{n},$$
where $I_n$ is the $n\times n$ identity matrix. These algebras,
first studied by Leavitt in \cite{Leavitt57} and \cite{Leavitt62},
provide universal examples of algebras without the invariant basis
number property: observe that right multiplication by the row
$(x_1,\dots ,x_n)$ gives an isomorphism from the free left
$L(1,n)$-module of rank one onto the free left $L(1,n)$-module of
rank $n$. They are algebraic analogues of the {\it Cuntz algebras}
$\mathcal O _n$, introduced independently by Cuntz in
\cite{Cuntz77}. The first author analyzed in \cite{Ara04} the
structure of the finitely presented modules over $L(1,n)$ in
connection with the structure of certain classes of finitely
presented modules over the free algebras $\cb \langle x_1,\dots
,x_n\rangle$ and $\cb \langle y_1,\dots ,y_n \rangle$. Both free
algebras embed in $L(1,n)$, and the abelian category $\mathcal S$ of
finitely presented left $L(1,n)$-modules of finite length is
equivalent to a quotient category of the abelian category of
finite-dimensional $k\langle y_1,\dots ,y_n\rangle$-modules by a
certain Serre subcategory, see \cite[Theorem 5.1]{Ara04}. Let
$\Sigma$ be the class of all the square matrices over $\cb \langle
x_1,\dots ,x_n \rangle$ that are sent to an invertible matrix by the
augmentation map. Then $\mathcal S$ is identified with the category
of finitely presented $\Sigma$-torsion modules in \cite[Theorem
6.2]{Ara04}, and this is used to give a formula for $K_1(Q_n)$,
where $Q_n=L(1,n)\Sigma^{-1}$ is the universal localization of
$L(1,n)$ with respect to $\Sigma$, which was shown in
\cite{AraGoodearlPardo02} to be a simple von Neumann regular ring.

The main purpose of this paper is to generalize these results to the
much wider context of path algebras. Our main guiding principle in
tackling this problem is the idea that free algebras are
prototypical examples of path algebras, and many results on free
algebras should admit suitable generalizations to this setting. For
each finite (or even row-finite) quiver $E$, there is a {\it Leavitt
path algebra} $L_\cb (E)$, described below, which plays a similar
role with respect to the usual path algebra $P_\cb (E)$ as $L(1,n)$
does with respect to the free algebra $\cb \langle x_1,\dots ,x_n
\rangle$. (Recall that $\cb \langle x_1,\dots ,x_n \rangle$ is the
path algebra of the quiver with one vertex and $n$ arrows.) The
Leavitt path algebras $L_\cb (E)$ were first introduced in
\cite{AbramsAranda05} and \cite{AraMorenoPardo07}, and have been
intensively studied by various authors since then. The {\it regular
algebra of} $E$, denoted by $Q_\cb (E)$, was constructed in
\cite{AraBrustenga07}, and is the natural generalization of the
algebra $Q_n$ described above; see below for the definition. It
follows from \cite[Theorem 4.2]{AraBrustenga07} that $K_0(Q_\cb
(E))\cong K_0(L_\cb (E))$ for every finite quiver $E$. We will
compute here (Theorem \ref{theor:K1QE}) all the higher $K$-theory
groups of $Q_\cb (E)$ in terms of the $K$-theory groups of $L_\cb
(E)$, recently computed in \cite{AraBrustengaCortinas09}, and the
$K$-theory of a certain abelian category $\Bla (\cam{E})$ of objects
of finite length. This is new even for the regular algebra $Q_n$ of
the classical Leavitt algebra $L(1,n)$, since only $K_1$ was
considered in \cite{Ara04}.

\vskip.3cm

Unless otherwise is stated all modules are left modules. In the
following, $\cb$ will denote a field and $E=(E^0,E^1,r,s)$ a finite
quiver (oriented graph) with $E^0=\{1,\dotsc,d\}$. Here $s(e)$ is
the {\em source vertex} of the arrow $e$, and $r(e)$ is the {\em
range vertex} of $e$. A {\em path} in $E$ is either an ordered
sequence of arrows $\alpha=e_1\dotsb e_n$ with $r(e_t)=s(e_{t+1})$
for $1\leqslant t<n$, or a path of length $0$ corresponding to a
vertex $i\in E^0$, which will be denoted by $p_i$. The paths $p_i$
are called trivial paths, and we have $r(p_i)=s(p_i)=i$. A
non-trivial path $\alpha=e_1\dotsb e_n$ has length $n$ and we define
$s(\alpha)=s(e_1)$ and $r(\alpha)=r(e_n)$. We will denote the length
of a path $\alpha$ by $|\alpha|$, the set of all paths of length $n$
by $E^n$ (for $n>1$), and the set of all paths by $E^*$.

Let us recall the construction of the Leavitt path algebra
$L(E)=L_\cb (E)$ and of the regular algebra $Q(E)=Q_{\cb}(E)$ of a
quiver $E$. These algebras fit into the following all-important
commutative diagram of injective algebra morphisms:

\begin{equation}
\label{maindiagram} \begin{CD}
\cb^d @>>> \cam{E} @>{\iota _{\Sigma}}>> \rat{E} @>>> \ser{E}\\
@VVV @V{\iota _{\Sigma_1}}VV @V{\iota _{\Sigma _1}}VV @V{\iota _{\Sigma _1}}VV\\
\cam{\ol{E}} @>\iota_{\Sigma_2}>> L(E) @>{\iota _{\Sigma}}>> Q(E)
@>>> U(E)
\end{CD}  \end{equation}
Here $P(E)$ is the path $\cb$-algebra of $E$, $\ol{E}$ denotes the
inverse quiver of $E$, that is, the quiver obtained by changing the
orientation of all the arrows in $E$, $\ser{E}$ is the algebra of
formal power series on $E$, and $\rat{E}$ is the algebra of rational
series, which is by definition the division closure of $P(E)$ in
$P((E))$ (which agrees with the rational closure, see
\cite[Observation 1.18]{AraBrustenga07}). The maps $\iota_{\Sigma}$
and $\iota _{\Sigma_i}$ indicate universal localizations with
respect to the sets $\Sigma$ and $\Sigma _i$ respectively. Here
$\Sigma$ is the set of all square matrices over $P(E)$ that are sent
to invertible matrices by the augmentation map $\epsilon \colon
P(E)\to \cb^{|E^0|}$, which coincides with the set of square
matrices over $\cam{E}$ which are invertible over $\ser{E}$
(\cite[Observation 1.19]{AraBrustenga07}). By \cite[Theorem
1.20]{AraBrustenga07}, the algebra $\rat{E}$ coincides with the
universal localization $P(E)\Sigma ^{-1}$. The set
$\Sigma_1=\{\mu_i\mid i\in E^0,\,s^{-1}(i)\neq \emptyset\}$ is the
set of morphisms between finitely generated projective left
$P(E)$-modules defined by
 \begin{align*}
  \mu_i\colon P(E)p_i&\longrightarrow \bigoplus_{j=1}^{n_i}P(E)p_{r(e^i_j)}\\
  r&\longmapsto\left(re^i_1,\dotsc,re^i_{n_i}\right)
 \end{align*}
for any $i\in E^0$ such that $s^{-1}(i)\neq\emptyset$, where
$s^{-1}(i)=\{ e_1^i,\dots ,e_{n_i}^i \}$. By a slight abuse of
notation, we use also $\mu _i$ to denote the corresponding maps
between finitely generated projective left $\rat{E}$-modules and
$\ser{E}$-modules respectively. The set $\Sigma_2=\{\nu_i\mid i\in
E^0,\,s^{-1}(i)\neq \emptyset\}$ is the set of morphisms between
finitely generated projective left $P(\ol{E})$-modules defined by
\begin{align*}
\nu_i\colon\bigoplus_{j=1}^{n_i} \cam{\overline{E}} p_{r(e^i_j)} &\longrightarrow \cam{\overline{E}}p_i\\
(r_1,\dotsc,r_{n_i})&\longmapsto\sum_{j=1}^{n_i}r_j
\overline{e}^i_j.
\end{align*}
for each $i\in E^0$ such that $s^{-1}(i)\neq\emptyset$.

The following relations hold in $Q(E)$:

(V)\, \, \, \,\hskip.35cm $p_vp_{v'}=\delta _{v,v'}p_v$  for all
$v,v'\in E^0$.

(E1) \, \, \hskip.3cm  $p_{s(e)}e=ep_{r(e)}=e$ for all $e\in E^1$.

(E2) \, \, \hskip.3cm $p_{r(e)}\ol{e}=\ol{e}p_{s(e)}=\ol{e}$ for all
$e\in E^1$.

(CK1)\hskip.5cm $\ol{e}e'=\delta _{e,e'}p_{r(e)}$ for all $e,e'\in
E^1$.

(CK2)\hskip.5cm  $p_v=\sum _{\{ e\in E^1\mid s(e)=v \}}e\ol{e}$ for
every $v\in E^0$ that emits edges.

\noindent The Leavitt path algebra
$L(E)=P(E)\Sigma_1^{-1}=P(\ol{E})\Sigma_2^{-1}$ is the algebra
generated by $\{p_v\mid v\in E^0\}\cup \{e,\ol{e}\mid e\in E ^1\}$
subject to the relations (V)--(CK2) above; see for instance
\cite{AbramsAranda05} and \cite{AraMorenoPardo07}. Relations (CK1)
and (CK2) are called the {\it Cuntz-Krieger relations}, see
\cite{CuntzKrieger80}. By \cite[Theorem 4.2]{AraBrustenga07}, the
algebra $Q(E)$ is a von Neumann regular hereditary ring and
$Q(E)=P(E)(\Sigma \cup \Sigma_1)^{-1}$.

A {\it sink} in $E$ is a vertex $i\in E^0$ such that
$s^{-1}(i)=\emptyset$, that is, $i$ does not emit any arrow. The set
of sinks of $E$ will be denoted by $\Si (E)$. With this terminology
we can summarize the results on the $K$-theory of the Leavitt
algebra $L_\cb (E)$, obtained in \cite{AraBrustengaCortinas09},  as
follows. Consider the adjacency matrix
$A_E=(a_{ij})\in\Z^{(E_0\times E_0)}$, $a_{ij}=\#\{$ arrows from $i$
to $j\}$. Write $N_E$ and $1$ for the matrices in $ \Z^{(E_0\times
E_0\setminus\Si(E))}$ which result from $A^t_E$ and from the
identity matrix after removing the columns corresponding to sinks.
Then there is a long exact sequence ($n\in \Z$)
\begin{equation*}
\begin{CD}
\cdots \to K_n(k)^{(E_0\setminus \Si (E))} @>{1-N_E}>>
K_n(k)^{(E_0)} @>>> K_n(L_k(E)) @>>> K_{n-1}(k)^{(E_0\setminus \Si
(E))}  .
\end{CD}
\end{equation*}
In particular
$$K_0(L_k(E))\cong \coker (1-N_E\colon \Z^{(E_0\setminus \Si
(E))}\longrightarrow \Z^{(E_0)}),$$ and
\begin{equation*}\begin{CD}K_1(L_k(E)) & \cong \coker (1-N_E\colon
(k^\times)^{(E_0\setminus \Si (E))}\longrightarrow (k^\times )^{(E_0)})\\
& \bigoplus \ker (1-N_E\colon \Z^{(E_0\setminus \Si
(E))}\longrightarrow \Z^{(E_0)}).\end{CD}\end{equation*} In Theorem
\ref{theor:K1QE}, we show that, for $i\ge 1$,
$$K_i(Q(E))\cong K_i(L(E))\bigoplus {\rm Bla}_{i-1}(P(E)),$$
where $\text{Bla}_*(\cam{E})$ is the $K$-theory of the abelian
category $\Bla (\cam{E})$ consisting of finitely generated
$\cam{E}$-modules $M$ such that ${\rm
Tor}_q^{P_\cb(E)}(\cb^{|E^0|},M)=0$ for all $q$. This category is
shown in Proposition \ref{typeLL} to be exactly the category of
finitely presented $\lev{E}$-modules of finite length without
nonzero projective submodules. Observe that, by the ``Devissage"
Theorem (\cite[5.3.24]{Rosenberg94}) and the results in the present
paper, the groups $\text{Bla}_i (\cam{E}))$ are the direct sum of
the $K_i$ groups of the endomorphism rings
$\text{End}_{\cam{\ol{E}}}(M)^{\text{op}}$, where $M$ ranges over
all the finite-dimensional non-projective simple
$\cam{\ol{E}}$-modules which are not isomorphic to one of the simple
modules $\coker (\nu_j)$ for $\nu_j\in \Sigma _2$.

The rest of the paper is organized as follows. As a preparation for
our main results, we develop in Sections 2 and 3 some results about
the structure of finitely presented modules over a path algebra.
This is done by extending to this context some of the tools
developed by Cohn to study firs. In particular we show in Theorem
\ref{tma:subproj} that every finitely related $\cam{E}$-module $L$
has a projective submodule $Q$ such that $L/Q$ is finite-dimensional
over $k$, generalizing a result of Lewin \cite{Lewin69} for the free
algebra. Section 4 establishes the important fact that $L(E)$ is
flat as a {\it right} $P(\ol{E})$-module, which will be often used
afterwards. We start our study of the module theory over Leavitt
path algebras in Section 5, obtaining in Proposition
\ref{prop:modfp} a description of the finitely presented
$\lev{E}$-modules of finite length as induced modules from
finite-dimensional $\cam{\ol{E}}$-modules. In Section 6, the abelian
categories $\fp{\lev{E}}$ and $\fpl{\lev{E}}$ of finitely presented,
and finitely presented $\lev{E}$-modules of finite length,
respectively, are shown to be equivalent to the quotient categories
of the corresponding categories of $\cam{\ol{E}}$-modules modulo the
Serre subcategory generated by the simple finite-dimensional
$\cam{\ol{E}}$-modules $\coker (\nu_j)$, for $\nu _j\in \Sigma _2$.
Finally we discuss the notion of Blanchfield modules in Section 7,
which we have adapted from \cite{RanickiSheiham06}, and we show that
the category of finitely generated Blanchfield $\cam{E}$-modules
agrees with various relevant categories. In particular it is the
category of torsion modules for both universal localizations
$\cam{E}\to \cam{E}\Sigma^{-1}$ and $\lev{E}\to \lev{E}\Sigma^{-1}$
(Proposition \ref{torsions}), and coincides with the category of
finitely presented $\lev{E}$-modules of finite length without
nonzero projective submodules (Proposition \ref{typeLL}). The
$K$-theory results described above are deduced then from the long
exact sequence of Neeman and Ranicki for stably flat universal
localizations
\cite{NeemanRanicki01},\cite{NeemanRanicki04},\cite{Neeman07}.

\section{Finitely presented modules over path algebras}
Let $\cb $ be a field and let $R=\cb\left<X\right>$ be the free
algebra in $n$ variables. Recall that given an $R$-module $M$ of
finite $\cb$-dimension we have the Lewin-Schreier formula relating
$\chi_R(M)$, the Euler characteristic, with the $\cb$-dimension of
$M$:
\[
 \chi_R(M)=(1-n)\dim_{\cb}(M)
\]
(see \cite[Theorem~4]{Lewin69} or \cite[Theorem~2.5.3]{Cohn85}).
Using a general result due to Bergman and Dicks
\cite{BergmanDicks78} we will see that a similar formula holds for
the path algebra.

To state the formula in our situation we will need a more general
context. Let $R$ be any ring. If an $R$-module $M$ has a finite
resolution by finitely generated projective modules,
\[
 0\longrightarrow P_n\longrightarrow\dotsb\longrightarrow P_0\longrightarrow M \longrightarrow 0,
\]
it is known that the element $\chi _R(M):=\sum(-1)^i[P_i]\in K_0(R)$
is an invariant of $M$ called its {\em Euler characteristic}.

Let $A$ be any ring. If $R$ is an $A$-ring, then it makes sense to compare $\chi_A(M)\in K_0(A)$ and $\chi_R(M)\in K_0(R)$ when both are defined. We have the following definition due to Bergman and Dicks:
\begin{defi}[{\cite[(64)]{BergmanDicks78}}]
 An $A$-ring $R$ will be called a {\it left Lewin-Schreier $A$-ring} if
\begin{enumerate}
 \item every left $R$-module $M$ which has a finite resolution by finitely generated projectives over
 $A$ also has such a resolution over $R$, and
 \item there exists a homomorphism $\lambda_{R}^{A}\colon K_0(A)\to K_0(R)$ such that, for such an $M$,
 $\chi_R(M)=\lambda_{R}^{A}\chi_A(M)$.
\end{enumerate}
\end{defi}

Let $R$ be an $A$-ring. We will denote by $\tau_{R}^{A}\colon
K_0(A)\to K_0(R)$ the homomorphism induced by the functor $R
\otimes_A -$.
\begin{prop}\label{prop:LewinSchreier}
Let $E$ be a finite quiver with $E^0=\{1,\dotsc,d\}$. Then $\cam{E}$
is a left Lewin-Schreier $\cb^d$-ring with
$\lambda_{\cam{E}}^{\cb^d}=(\id-A_E^t)\tau_{\cam{E}}^{\cb^d}$.
\end{prop}
\begin{proof}
% We will proceed by induction on $n=E^1$. Let $E$ be a quiver with a single edge $e$ (and $d$ vertices); we write $R=\cam{E}$ and $A=\cb^d$. Let $M\subseteq R$ be the $A$-bimodule generated by $e$. It is easy to check that the path algebra of a quiver is isomorphic to the tensor $A$-ring associated to the bimodule generated by the edges (see \cite[Proposition~III.1.3]{AuslanderReitenSmalo95}). Therefore, by \cite[(63)]{BergmanDicks78} we get the following exact sequence:
% \begin{equation}\label{eq:es}
% 0\longrightarrow R\otimes_AM\otimes_AR\longrightarrow R\otimes_AR\longrightarrow R\longrightarrow 0.
% \end{equation}
% Let $N$ be a right $R$-module finitely generated as $A$-module. Applying the functor $N\otimes_R-$ to the exact sequence~\eqref{eq:es} we get a resolution of $N$ by finitely generated projective $R$-modules
% \[
% 0\longrightarrow N\otimes_AM\otimes_AR\longrightarrow N\otimes_AR\longrightarrow N\longrightarrow 0.
% \]
% So $R$ satisfies the first condition in the definition.
%
% As $A$-module, $N$ is isomorphic to $(p_1A)^{\alpha_1}\oplus\dotsb\oplus(p_dA)^{\alpha_d}$ for some $\alpha_1,\dotsc,\alpha_d\in\N$. For $i=1,\dotsc,d$, we have the following isomorphisms of right $R$-modules
% \[
% p_iA\otimes_AR\cong p_iR,\qquad p_iA\otimes_AM\otimes_AR=\begin{cases}0\quad&\text{if $i\neq s(e)$}\\p_{r(e)}R\quad&\text{if $i=s(e)$.}\end{cases}
% \]
% So, we get
% \begin{multline*}
% \chi_R(N)=[N\otimes_A R]-[N\otimes_A M\otimes_AR]=\sum_{i=1}^d\alpha_i[p_iR]\\
% \end{multline*}
We write $R=\cam{E}$ and $A=\cb^d$. Let $N\subseteq R$ be the
$A$-bimodule generated by the edges. It is easy to check that the
path algebra of a quiver is isomorphic to the tensor $A$-ring
associated to the bimodule generated by the edges (see
\cite[Proposition~III.1.3]{AuslanderReitenSmalo95}). Therefore, by
\cite[(63)]{BergmanDicks78} we get the following exact sequence:
\begin{equation}\label{eq:es}
0\longrightarrow R\otimes_AN\otimes_AR\longrightarrow
R\otimes_AR\longrightarrow R\longrightarrow 0.
\end{equation}
Let $M$ be a left $R$-module finitely generated as $A$-module.
Applying the functor $- \otimes_R M$ to the exact
sequence~\eqref{eq:es} we get a resolution of $M$ by finitely
generated projective left $R$-modules
\[
0\longrightarrow R\otimes_AN\otimes_AM\longrightarrow
R\otimes_AM\longrightarrow M\longrightarrow 0,
\]
and so, $R$ satisfies the first condition in the definition.

As an $A$-module, $M$ is isomorphic to
$(Ap_1)^{\alpha_1}\oplus\dotsb\oplus(Ap_d)^{\alpha_d}$ for some
$\alpha_1,\dotsc,\alpha_d\in\N$. We put $A_E=(a_{ij})$. We have the
following isomorphisms of left $R$-modules
\[
R\otimes_AAp_i\cong Rp_i,\qquad
R\otimes_AN\otimes_AAp_i\cong\bigoplus_{j=1}^d(Rp_j)^{a_{ji}}.
\]
So, we get
\[
\chi_R(M)=[R\otimes_A M]-[R\otimes_A N\otimes_A
M]=\sum_{i=1}^d\alpha_i[Rp_i]-
\sum_{i=1}^d\sum_{j=1}^da_{ji}\alpha_j[Rp_i]= (\id-A_E^t)\tau_{R}^{A}\chi_A(M)\\
\]
as wanted.
\end{proof}

\section{The weak algorithm for path algebras}

The path algebra can be profitably thought of as a generalization of the free algebra and,
quite often, properties of the latter admit a generalization to the former.
In this section we generalize Cohn's weak algorithm (see \cite[Chapter~2]{Cohn85}) to
the context of path algebras and prove several of its basic properties.
%Beware that, although our definitions are slightly different from Cohn's ones, we will keep the same names for analogous concepts.
The main result in this section is Theorem~\ref{tma:subproj} which is a version of
Lewin's Theorem (see \cite[Theorem~2]{Lewin69}) for path algebras.

Let $R$ be a non-zero ring. Recall that a {\em filtration} on $R$ is
given by a map $\nu\colon R\to \N\cup\{-\infty\}$ with the following
properties:
\begin{enumerate}
\item $\nu(r)\geqslant 0$ for all $r\neq 0,\,\nu(0)=-\infty$,
\item $\nu(r-s)\leqslant\max\{\nu(r),\nu(s)\}$,
\item $\nu(rs)\leqslant\nu(r)+\nu(s)$,\label{enu1:item3}
\item $\nu(1)=0$.
\end{enumerate}
If equality holds in \eqref{enu1:item3}, we have a degree function. Even in the general case we shall call $\nu(r)$ the {\em degree} of $r$. It is easy to see that the path algebra $\cam{E}$ is a filtered ring with respect to the degree. A filtration is also determined by the additive subgroups $R_h$ given by the elements of degree at most $h$.

% Given a filtration, let us write $R_h$ for the set of elements of degree at most $h$; then the $R_h$ are subgroups of the additive group of $R$ such that
% \begin{enumerate}[(i)]
% \item $0=R_{-\infty}\subset R_0\subseteq R_1\subseteq\dotsc,$ \label{enu2:item1}
% \item $\cup R_h=R,$
% \item $R_iR_j\subseteq R_{i+j},$
% \item $1\in R_0.$\label{enu2:item4}
% \end{enumerate}
% Conversely, every series of subgroups $R_h$ of the additive group of $R$ satisfying \eqref{enu2:item1}--\eqref{enu2:item4} leads to a filtration $\nu$, given by $\nu(r)=\min\{h\mid r\in R_h\}$.

Let $R$ be a ring with a filtration $\nu$. Given an $R$-module $M$ a
{\em filtration} on $M$ is given by a map $\mu\colon M\to
\N\cup\{-\infty\}$ such that
\begin{enumerate}
\item $\mu(m)\geqslant 0$ for all $m\neq 0,\,\mu(0)=-\infty$,
\item $\mu(m-n)\leqslant\max\{\mu(m),\mu(n)\}$,
\item $\mu(mr)\leqslant\mu(m)+\nu(r)$.
\end{enumerate}
Like in the ring case, a filtration on $M$ is also determined by the
additive subgroups $M_h$ given by the elements of degree at most
$h$.

% let $M_h$ denote the set of elements of degree at most $h$; these are additive subgroups of $M$ satisfying:
% \begin{enumerate}[(i)]
% \item $0=M_{-\infty}\subseteq M_0\subseteq M_1\subseteq\dotsc,$
% \item $\cup M_h=M,$
% \item $M_iR_j\subseteq M_{i+j}.$
% \end{enumerate}
% As before, given subgroups $M_h$ of $M$ satisfying the latter conditions leads to a filtration $\mu$ given by $\mu(m)=\min\{h\mid m\in M_h\}$.

The following definition is useful to generalize Cohn's concept of $\mu$-independence to the context of path algebras.
\begin{defi}
Let $(R,\nu)$ be a filtered ring. A {\em set of vertices in} $R$ is
a finite set $P$ of zero-degree, pairwise orthogonal idempotents in
$R$ such that $1=\sum_{p\in P}p$. We also say that $R$ has a {\em
vertex-type decomposition given by} $P$.
\end{defi}
\begin{exems}
\begin{enumerate}[1.]
 \item Any filtered ring has a trivial vertex-type decomposition given by $P=\{1\}$.
 \item The path algebra of a finite quiver $E$ has a vertex-type decomposition given by the vertices $P=\{p_i\mid i\in E^0\}$. This is the example to bear in mind.
 \item Mixed path algebras as defined in \cite{AraBrustenga08} have also a vertex-type decomposition given by the vertices.
\end{enumerate}
\end{exems}

In the following definitions and results $R$ will denote a ring with
a filtration $\nu$, $P=\{p_1,\dotsc,p_d\}$ will be a set of vertices
in $R$ and $M$ will be an $R$-module with a filtration $\mu$.

\begin{defi}
We say that the family $(m_i)_{i\in I}\in\prod_{i\in I}p_{n_i}M$ is
{\em left $P$-$\mu$-dependent} provided that exists a family
$(r_i)_{i\in I}\in \bigoplus_{i\in I} Rp_{n_i}$ such that
\[
\mu\left(\sum_{i\in I} r_im_i\right)<\max_{i\in
I}\{\nu(r_i)+\mu(m_i)\}
\]
or if some $m_i=0$. Otherwise the family $(m_i)_{i\in I}$ is said to
be {\em left $P$-$\mu$-independent}.
\end{defi}

When $P=\{1\}$ (and $M=R$) we recover Cohn's definitions of
left $\mu$-dependent and left $\mu$-independent family (see \cite[Pag.~95]{Cohn85}).
Recall that in Cohn's setting a left $\mu$-independent family generates a free module
(because it is also a left linearly independent family). In the general case, the point
is the fact that a left $P$-$\mu$-independent family generates a projective module: %({\em cf} Lemma~\ref{lem:quocients} below).

\begin{prop}\label{prop:indeproj}
In the above situation, let $(m_i)_{i\in I}\in\prod_{i\in
I}p_{n_i}M$ be a $P$-$\mu$-independent family. Then the submodule
$\sum_{i\in I}Rm_i$ is projective.
\end{prop}
\begin{proof}
Indeed, by the $P$-$\mu$-independence of the family, the epimorphism
 \begin{align*}
  \bigoplus_{i\in I}Rp_{n_i}&\longrightarrow \sum_{i\in  I}Rm_i\subseteq M\\
  (r_i)_{i\in I}&\longmapsto \sum_{i\in I}r_im_i
 \end{align*}
 is an isomorphism.
\end{proof}

\begin{defi}
An element $m\in M$ is said to be {\em left $P$-$\mu$-dependent} on
a family $(m_i)_{i\in I}\in\prod_{i\in I} p_{n_i}M$ if either $m=0$
or there exists a family $(r_i)_{i\in I}\in\bigoplus_{i\in I}
Rp_{n_i}$ such that
\[
\mu\left(m-\sum_{i\in I} r_im_i\right)<\mu(m)\quad\text{and}\quad
\forall i\in I,\ \nu(r_i)+\mu(m_i)\leqslant \mu(m).
\]
In the contrary case $m$ is said to be {\em left
$P$-$\mu$-independent} of $(m_i)_{i\in I}$.
\end{defi}

We will also need the definition of left $P$-$\mu$-dependence of an
element on a general set:
\begin{defi}\label{def:genemudep}
An element $m\in M$ is said to be {\em left $P$-$\mu$-dependent} on
a set $S\subseteq M$ provided that there exists a family
$(m_i)_{i\in I}\in\prod_{i\in I} p_{n_i}S$ such that $m$ is left
$P$-$\mu$-dependent on it. Otherwise $m$ is said to be {\em left
$P$-$\mu$-independent} of $S$.
\end{defi}

Now, we can generalize the weak algorithm to our framework:
\begin{defi}
We say that $M$ satisfies the {\em weak algorithm} relative to $\mu$
and $P$ if in every finite left $P$-$\mu$-dependent family
$(m_i)_{i=1,\dotsc,\ell}\in\prod_{i=1}^\ell p_{n_i}M$ where
\[
\mu(m_1)\leqslant\dotsb\leqslant\mu(m_\ell),
\]
some $m_i$ is left $P$-$\mu$-dependent on $m_1,\dotsc,m_{i-1}$.
\end{defi}

Applying these definitions to the regular module $M={_R}R$ with the
filtration $\mu=\nu$ we also have these concepts defined for the
filtered ring $(R,\nu)$.

%Given an expression $\sum_{i\in I}m_ir_i\in M$ with $m_i\in Mp_{n_i}$ and $r_i\in p_{n_i}R$ we will refer to $\max_{i}\{\mu(m_i)+\nu(r_i)\}$ as its {\em formal degree}. We remark that the definition of $P$-$\mu$-independence of a family states that the degree of elements represented by certain expressions should equal the formal degree of these expressions. For a general expression $\sum_{i\in I}m_ir_i$ we will take its formal degree to be the formal degree of $\sum_{i\in I}\sum_{j=1}^d(m_ip_j)(p_jr_i)$.
Given an expression $\sum_{i\in I}r_im_i\in M$ with $m_i\in M$ and
$r_i\in R$ we will refer to $\max_{i}\{\nu(r_i)+\mu(m_i)\}$ as its
{\em formal degree}. We remark that the definition of
$P$-$\mu$-independence of a family states that the degree of
elements represented by certain expressions should equal the formal
degree of these expressions.

The previous definitions are motivated by the fact that any free
module over the path algebra satisfies the weak algorithm relative
to a suitable degree, as we show in our next result. This will be
improved in Theorem \ref{thm:projectives}, where it is shown that
the $\cam{E}$-modules satisfying the weak algorithm relative to some
filtration are precisely the projective $\cam{E}$-modules.
\begin{prop}\label{prop:WA}
Let $E$ be a finite quiver with $E^0=\{1,\dotsc,d\}$. Let $M$ be a
free $\cam{E}$-module freely generated by $\mathcal{B}$ and consider
a map $\mu\colon\mathcal{B}\to\N$. If we extend $\mu$ to $M$ as the
formal degree, then $(M,\mu)$ is a filtered module and satisfies the
weak algorithm relative to $\mu$ and $P=\{p_1,\dotsc,p_d\}$, the set
of vertices given by the vertices of $E$.
\end{prop}
\begin{proof}
First of all, since elements in $M$ have a unique expression as
$\cam{E}$-linear combination of elements in $\mathcal{B}$, the
formal degree gives a well defined filtration on $M$. Now we will
prove that $M$ satisfies the weak algorithm relative to $\mu$ and
$P$. Let $(m_i)_{i=1,\dotsc,\ell}\in\prod_{i=1}^\ell
(p_{n_i}M\setminus\{0\})$ be a left $P$-$\mu$-dependent family such
that $\mu(m_1)\leqslant\dotsb\leqslant\mu(m_\ell)$. There exists an
element $(r_i)_{i=1,\dotsc,\ell}\in \bigoplus_{i=1}^\ell
\cam{E}p_{n_i}$ such that
\begin{equation}\label{eq:mudep}
\mu\left(\sum_{i=1}^\ell
r_im_i\right)<t=\max_i\{\nu(r_i)+\mu(m_i)\}.
\end{equation}
By omitting some terms if necessary we may assume that, for all $i$,
$\nu(r_i)+\mu(m_i)=t$ and hence
$\nu(r_\ell)\leqslant\dotsb\leqslant\nu(r_1)$.

Since $\mathcal{B}$ is a basis for $M$, every $m_i$ has a unique
expression $m_i=\sum_{b\in\mathcal{B}}r^i_bb$. Moreover, from
$p_{n_i}m_i=m_i$ we get that $p_{n_i}r_b^i=r_b^i$. Therefore,
\begin{equation*}%\label{eq:filform}
\mu\left(\sum_{i=1}^\ell r_im_i\right)=\mu\left(\sum_{i=1}^\ell
r_i\left(\sum_{b\in\mathcal{B}}r_b^ib\right)\right)=
\mu\left(\sum_{b\in\mathcal{B}}\left(\sum_{i=1}^\ell r_ir_b^i\right)b\right).%\max_{b\in\mathcal{B}}\left\{\mu(b)+\nu\left(\sum_{i=1}^\ell r_b^ir_i\right)\right\}.
\end{equation*}

Let $\gamma\in\supp{r_\ell}$ (the support of $r_\ell$) be a path of maximal length, say $t_0$.
Now, given $r,s\in \cam{E}$, we have that
\begin{equation}\label{eq:cong}
\delta_\gamma(sr)\equiv \delta_\gamma(s)r\pmod {\cam{E}_{\nu(r)-1}},
\end{equation}
where $\delta_\gamma $ is the right transduction corresponding to
$\gamma$, that is, $\delta_\gamma (\gamma \tau ')=\tau '$ and
$\delta (\tau)=0$ if $\tau $ does not start with $\gamma$;  see
\cite[Section 1]{AraBrustenga07}. This is clear if $s$ is a monomial
of length at least $t_0$; in fact we then have equality. If $s$ is a
monomial of length less than $t_0$, the right-hand side of
\eqref{eq:cong} is zero, and so it holds as a congruence. The
general case follows by linearity.

Now, for all $i$ and all $b$, the element $\delta_\gamma (r_i)r^i_b$
differs from $\delta_\gamma (r_ir_b^i)$ by a term of degree less
than $\nu(r_b^i)$. Therefore, we have
\begin{equation*}
\nu\left(\sum_{i=1}^\ell\left(\delta_\gamma(r_i)r_b^i-
\delta_\gamma(r_ir_b^i)
\right)\right)\leqslant\max_i\{\nu(\delta_\gamma(r_i)r_b^i -
\delta_\gamma(r_ir_b^i))\}<\max_i\{\nu(r_b^i)\}.
\end{equation*}
From this inequality, we get
\begin{multline}\label{eq:dif}
\mu\left(\sum_{b\in\mathcal{B}}\left(\sum_{i=1}^\ell \delta_\gamma(r_i)r_b^i\right)b-
\sum_{b\in\mathcal{B}}\delta_\gamma \left(\sum_{i=1}^\ell r_ir_b^i\right)b\right)=\\
\qquad\shoveleft{=\mu\left(\sum_{b\in\mathcal{B}}\left(\sum_{i=1}^\ell\left(\delta_\gamma(r_i)r_b^i
-\delta_\gamma\left(r_ir_b^i\right)\right)\right)b\right)}\\
\qquad\shoveleft{=\max_{b\in\mathcal{B}}\left\{\mu(b)+\nu\left(\sum_{i=1}^\ell\left(\delta_\gamma(r_i)r_b^i-
\delta_\gamma\left(r_ir_b^i\right)\right)\right)\right\}}\\
\qquad\shoveleft{<\max_{b\in\mathcal{B}}\left\{\mu(b)+\max_i\{\nu(r_b^i)\}\right\}}\\
\qquad\shoveleft{=\max_{b\in\mathcal{B}}\left\{\max_i\{\mu(r_b^ib)\}\right\}
=\max_i\left\{\max_{b\in\mathcal{B}}\left\{\mu(r_b^ib)\right\}\right\}}\\
\qquad\shoveleft{=\max_i\left\{\mu\left(\sum_{b\in\mathcal{B}}r_b^ib\right)\right\}=
\max_i\{\mu(m_i)\}=\mu(m_\ell).}\\
\end{multline}
On the other hand, we have
\begin{align}
\mu\left(\sum_{b\in\mathcal{B}} \delta_\gamma\left(\sum_{i=1}^\ell
r_ir_b^i \right)b\right)&=
\max_{b\in\mathcal{B}}\left\{\mu(b)+\nu\left(\delta_\gamma\left(\sum_{i=1}^\ell
r_ir_b^i\right)\right)
\right\}\label{eq:mi}\\
\notag&\leqslant\max_{b\in\mathcal{B}}\left\{\mu(b)+\nu\left(\sum_{i=1}^\ell r_ir_b^i\right)\right\}-t_0\\
\notag&=\mu\left(\sum_{i=1}^\ell
r_im_i\right)-t_0\\
\notag&<t-t_0=\mu(m_\ell).
\end{align}
Hence, by \eqref{eq:dif} and \eqref{eq:mi} we get that
\[
\mu\left(\sum_{i=1}^\ell
\delta_\gamma(r_i)m_i\right)=\mu\left(\sum_{b\in\mathcal{B}}\left(\sum_{i=1}^\ell
\delta_\gamma(r_i)r_b^i\right) b \right)<\mu(m_\ell)
\]
and, since $\delta_\gamma(r_\ell)\in \cb^\times p_{n_{\ell}}$ we
deduce that $m_\ell$ is left $P$-$\mu$-dependent on
$m_1,\dotsc,m_{\ell-1}$ as wanted.
\end{proof}

In particular, the path algebra $\cam{E}$ satisfies the weak
algorithm relative to the degree and the obvious set of vertices. It
is straightforward to see that the weak algorithm is inherited by
submodules:
\begin{lem}\label{lem:restric}
Let $(R,\nu)$ be a filtered ring with a set of vertices $P$ and let
$(M,\mu)$ be a filtered right $R$-module satisfying the weak
algorithm relative to $\mu$ and $P$. Then every submodule
$N\subseteq M$ satisfies the weak algorithm relative to $\mu_{|N}$
and $P$.
\end{lem}
% \begin{proof}
%  Let $(m_i)_{i=1,\dotsc,\ell}\in\prod_{i=1}^\ell (Np_{n_i}\setminus\{0\})$ be a family right $P$-$\mu'$-dependent such that $\mu'(m_1)\leqslant\dotsb\leqslant\mu'(m_\ell)$. Therefore, the same is true for $\mu$ and, by the weak algorithm, some $m_k$ is $P$-$\mu$-dependent on $m_1,\dotsc,m_{k-1}$ and hence $P$-$\mu'$-dependent.
% \end{proof}

We have the following restriction for rings with weak algorithm:
\begin{prop}
Let $(R,\nu)$ be a filtered ring with a set of vertices $P$. If $R$
satisfies the weak algorithm relative to $\nu$ and $P$ then $R_0$ is
a semisimple ring.
\end{prop}
\begin{proof}
The set $R_0=\{r\in R\mid \nu(r)\leqslant 0\}$ is clearly a subring
of $R$. We have a finite decomposition $R_0=\bigoplus_{p\in P}R_0p$
into left ideals and we just need to check that these are simple
ideals. Fix some $p\in P$, since $\nu(p)=0$ we see that $R_0p$ is a
non-zero left ideal. Let $r\neq 0$ be in $R_0p$ and pick $q\in P$
such that $qr\neq 0$. Now the pair $(qr,p)$ is left
$P$-$\nu$-dependent and, by the weak algorithm, $p$ is left
$P$-$\nu$-dependent on $qr$, i.e. there exists $s\in R_0q$ such that
$\nu(p-sqr)<\nu(p)=0$. Thus $sqr=p$ and $R_0p$ is simple.
\end{proof}

\begin{defi}
Let $(R,\nu)$ be a filtered ring with a set of vertices $P$ and let
$(M,\mu)$ be a filtered $R$-module. A subset $\mathcal{B}$ of
$\cup_{p\in P}pM$ will be called a {\em weak $P$-$\mu$-basis} for
$M$ provided that
\begin{enumerate}[(i)]
\item Every element in $M$ is left $P$-$\mu$-dependent on $\mathcal{B}$.
\item No element of $\mathcal{B}$ is left $P$-$\mu$-dependent on the rest of $\mathcal{B}$.\label{enu3:item2}
\end{enumerate}
\end{defi}
It is easily seen, using the well-ordering of the range of $\mu$,
that a weak $P$-$\mu$-basis of $M$ generates $M$ as an $R$-module;
but in general it need be neither $P$-$\mu$-independent nor a
minimal generating set. However if $M$ satisfies the weak algorithm
relative to $\mu$ and $P$ then every weak $P$-$\mu$-basis of $M$ is
left $P$-$\mu$-independent by condition \eqref{enu3:item2} and
hence, by Proposition~\ref{prop:indeproj}, the module $M$ is
projective.

%Applying the previous definitions to the regular module $M=R_R$ with the filtration $\mu=\nu$ we also have this concepts defined for the ring $R$.

The remaining results in this section work in a more general setting
but we will state them only for the path algebra, which is the case
that we are interested in. From now on $E$ will be a finite quiver
with $E^0=\{1,\dotsc,d\}$, $\nu$ will denote the usual degree in the
path algebra and $P=\{p_1,\dotsc,p_d\}$ will be the natural set of
vertices of the path algebra. We can assure existence of weak
$P$-$\mu$-basis for filtered $\cam{E}$-modules:

\begin{prop}\label{prop:Ebases}
Let $(M,\mu)$ be a filtered $\cam{E}$-module. Then there exist sets
$\mathcal{B}^i_h\subseteq p_iM_h\setminus M_{h-1}$, for all
$i=1,\dotsc,d$ and $h\in \N$, such that
$\mathcal{B}=\cup_{i,h}\mathcal{B}^i_h$ is a weak $P$-$\mu$-basis
for $M$. Moreover, the cardinality of $\mathcal{B}^i_h$ does not
depend on the weak $P$-$\mu$-basis.
\end{prop}

\begin{proof}
The additive subgroup $M_h=\{m\in M\mid \mu(m)\leqslant h\}$ has an
structure of $\cb^d$-module induced by the inclusion $\cb^d\subseteq
\cam{E}$. For $h>0$ we denote by $M_h'$ the set of elements in $M_h$
left $P$-$\mu$-dependent on the set $M_{h-1}$ and put $M_0'=\{0\}$.
Observe that $M_h'$ is also a $\cb^d$-module. Indeed, it is clear
that $M_h'$ is closed under left product by elements in $\cb^d$;
closure with respect to the sum is clear if it has degree $h$ and,
otherwise it belongs to $M_{h-1}$. So, we may consider the
$\cb^d$-module $M_h/M_h'$ and the set $p_i(M_h/M_h')$ is a
$\cb$-vector space. Now, for every $h\geqslant 0$ and $i=1,\dotsc,d$
we pick $\mathcal{B}^i_h\subseteq M_h$ a set of representatives for
a $\cb$-basis of $p_i(M_h/M_h')$ such that $\mathcal{B}_h^i\subseteq
p_iM_h$. We write $\mathcal{B}=\cup_{i,h}\mathcal{B}^i_h$.

We will show that $\mathcal{B}$ is a weak $P$-$\mu$-basis for $M$.
By induction on $h$ every element in $M_h$ is left
$P$-$\mu$-dependent on $\mathcal{B}$. Indeed, for $h=0$ this holds
by construction. Assume that the statement is true for $h\geqslant
0$. By construction, every element in $M_{h+1}$ differs in some
element in $M'_{h+1}$ from a $\cb^d$-linear combination of elements
in $\mathcal{B}$ (of degree $h+1$). Every element in $M_{h+1}'$ is
$P$-$\mu$-dependent on $M_h$ and every element in $M_h$ is
$P$-$\mu$-dependent on $\mathcal{B}$. Therefore every element in
$M_{h+1}$ is $P$-$\mu$-dependent on $\mathcal{B}$. Moreover, since
$M=\cup_h M_h$, every element in $M$ is $P$-$\mu$-dependent on
$\mathcal{B}$.

Suppose that there is $b\in\mathcal{B}$ left $P$-$\mu$-dependent on
$\mathcal{B}\setminus\{b\}$. We write $h=\mu(b)$ and let $p_j\in P$
be such that $p_jb=b$. By construction $b\neq 0$, and hence there
exist $(b_i)_{i\in I}\in\prod_{i\in
I}(p_{n_i}\mathcal{B}\setminus\{b\})$ and $(r_i)_{i\in
I}\in\bigoplus_{i\in I}Rp_{n_i}$ such that
\[
\mu\left(b-\sum_{i\in I}r_ib_i\right)<h\quad\text{and}\quad\forall
i\in I,\,\nu (r_i)+\mu(b_i)\leqslant h.
\]
Moreover, we can assume that, for all $i$, $p_jr_i=r_i$. For all $i$
such that $r_i\neq 0$ we have $\mu(b_i)\leqslant h$ and, if
$\mu(b_i)=h$ then $\nu(r_i)=0$, and so $p_{n_i}=p_j$; therefore $b$
differs in an element in $M_h'$ from a $\cb$-linear combination of
elements in $\mathcal{B}^{j}_h$. This contradicts the fact that
classes of elements in $\mathcal{B}^{j}_h$ are linearly independent
elements in $p_j(M_h/M_h')$. Thus, we get that $\mathcal{B}$ is a
weak $P$-$\mu$-basis for $M$.

On the other hand, given a weak $P$-$\mu$-basis $\mathcal{C}$ for
$M$ it is clear that classes modulo $M_h'$ of elements in the set
$\{c\in\mathcal{C}\mid p_ic=c,\,\mu(c)=h\}$ give a $\cb$-basis of
the $\cb$-vector space $p_i(M_h/M_h')$; hence, its cardinality does
not depend on the weak $P$-$\mu$-basis.
\end{proof}

Now we can characterize projective $\cam{E}$-modules using the weak algorithm:
\begin{theor}
\label{thm:projectives}
 A $\cam{E}$-module $M$ is projective if and only if $M$ satisfies the weak algorithm relative to a suitable filtration.
\end{theor}
\begin{proof}
Let $M$ be a projective $\cam{E}$-module. Then $M$ is a submodule of
some free $\cam{E}$-module, say $F$. By Proposition~\ref{prop:WA},
the free module $F$ satisfies the weak algorithm relative to some
filtration $\mu$ (and $P$). Therefore, by Lemma~\ref{lem:restric},
the module $M$ satisfies the weak algorithm relative to the
restriction $\mu_{|M}$.

Let $(M,\mu)$ be a filtered module satisfying the weak algorithm relative to $\mu$ and $P$.
By Proposition~\ref{prop:Ebases}, the module $M$ has a weak $P$-$\mu$-basis, which is
$P$-$\mu$-independent due to the weak algorithm. Hence, by Proposition~\ref{prop:indeproj},
the module $M$ is projective.
\end{proof}

Let $R$ be a ring and $M$ an $R$-module. Recall that $M$ is {\em finitely related} provided
that there is an exact sequence of $R$-modules
\[
 0\longrightarrow L \longrightarrow F\longrightarrow M\longrightarrow 0,
\]
where $F$ is a free module and $L$ is finitely generated.

The following result generalizes a Theorem by Lewin \cite[Theorem~2]{Lewin69}.
The idea of the proof lies on an unpublished demonstration of Lewin's result due to
Warren Dicks \cite{Dicks02}. We gratefully acknowledge him for providing it to us.
\begin{theor}\label{tma:subproj}
Let $L$ be a finitely related $\cam{E}$-module. Then $L$ contains a projective module $Q$ such that $L/Q$ has finite $\cb$-dimension.
\end{theor}
\begin{proof}
Let
\[
0\longrightarrow N\longrightarrow M\xrightarrow{\ \varphi\;} L\longrightarrow 0
\]
be a presentation for $L$, where $M$ is free on a subset
$\mathcal{E}$, say, and $N$ is a finitely generated submodule of
$M$. Moreover, since  $\cam{E}$ is a hereditary ring, $N$ is a
projective module. It is well-know (see e.g.
\cite[Proposition~1.2]{AraBrustenga07}) that $N$ is isomorphic to a
direct sum of copies of the modules $\cam{E}p_i$; hence, there
exists $(f_1,\dotsc,f_m)\in \prod_{i=1}^m p_{n_i}M$ such that
$\cam{E}f_i\cong \cam{E}p_{n_i}$ and
\[
 N=\bigoplus_{i=1}^m\cam{E} f_i\cong \bigoplus_{i=1}^m \cam{E}p_{n_i}.
\]  We write $\mathcal{F}=\{f_1,\dotsc,f_m\}$.

Elements in $\mathcal{F}$ are $\cam{E}$-linear combinations of elements
in $\mathcal{E}$. Consider a finite subset $\mathcal{E'}\subseteq\mathcal{E}$
such that expressions of elements in $\mathcal{F}$ only involve elements
in $\mathcal{E}'$. Now we define $\mu(\mathcal{E}')=1$ and extend $\mu$
to $\mathcal{F}$ as the formal degree determined by $\mu$ and $\nu$,
the degree in $\cam{E}$. We write $n=\max\{\mu(f)\mid f\in \mathcal{F}\}$,
define $\mu(\mathcal{E}\setminus\mathcal{E'})=n+1$ and extend $\mu$ to $M$
as the formal degree. By Proposition~\ref{prop:WA} we get
that $(M,\mu)$ satisfies the weak algorithm with respect to $\mu$ and $P$.

From Lemma~\ref{lem:restric}, $N$ also satisfies the weak algorithm
with respect to $\mu'=\mu_{|M}$ and, by
Proposition~\ref{prop:Ebases}, $N$ has a weak $P$-$\mu'$-basis, say
$\mathcal{F}'$. Therefore, $\mathcal{F}'$ is left
$P$-$\mu'$-independent. Moreover, since $N$ is finitely generated
and (by definition of $\mu'$) $P$-$\mu'$-dependent on $N_n$,
$\mathcal{F}'$ is finite and contained in $N_n$.

Now, we will construct a $P$-$\mu$-independent family in $M$ in such
a way that it gives rise to a projective submodule in $L$. We have the filtration $\mu''$ on $L$ determined by setting $L_h=(M_h+N)/N$ (viewing $L$ as $M/N$). %We define $L_{n+1}':=L_n$, let $\overline{\mathcal{B}}_{n+1}^i\subseteq L_{n+1}p_i$ be a set of representatives for a $\cb$-basis of $(L_{n+1}/L'_{n+1})p_i$ and write  $\overline{\mathcal{B}}_{n+1}=\cup_{i=1}^d\overline{\mathcal{B}}^i_{n+1}$. Now, for $t>n+1$ we define inductively the $\cb^d$-modules $L_t'$ and the sets $\overline{\mathcal{B}}_t$ in the following way: $L_t'$ is given by elements $p\in L_t$ such that $p$ is either $\mu''$-dependent on $\cup_{t>h>n}\overline{\mathcal{B}}_h$ or $\mu''(p)<t$ and $\overline{\mathcal{B}}_t$ is equal to $\cup_{i=1}^d\overline{\mathcal{B}}^i_t$ where $\overline{\mathcal{B}}^i_t\subseteq L_tp_i$ is a set of representatives for a $\cb$-basis of $(L_t/L_t')p_i$. For every $i\in E^0$ and $t>n$ \dots
Let $L_h'$ denote the set of elements of $L_h$ which are
$\mu''$-dependent on $L_{h-1}$. For $t>n$ and $i\in E^0$, let
$\mathcal{B}^i_t$ be a subset of $p_iM_t$ whose image is a
$\cb$-basis of $p_i(L_t/L_t')$. Write
$\mathcal{B}^i=\cup_{t>n}\mathcal{B}^i_t$,
$\mathcal{B}_t=\cup_{i=1}^d\mathcal{B}^i_t$ and
$\mathcal{B}=\cup_{i=1}^d\mathcal{B}^i=\cup_{t>n}\mathcal{B}_t$.
Consider the submodule
$Q=\sum_{i=1}^d\sum_{b\in\mathcal{B}^i}\cam{E}\varphi(b)\subseteq
L$. The $\cam{E}$-module epimorphism defined as follows
\begin{align*}
  \bigoplus_{i=1}^d\bigoplus_{b\in\mathcal B^i}\cam{E}p_i &\longrightarrow Q\\
  (r_b^i)_{i,b}&\longmapsto
  \sum_{i=1}^d\sum_{b\in\mathcal{B}^i}r_b^i\varphi(b)
\end{align*}
is an isomorphism. Indeed, suppose not, then there exist elements
$r^i_b\in \cam{E}p_i$ not all zero such that
$\sum_{i=1}^d\sum_{b\in\mathcal{B}^i}r^i_bb\in N$. Therefore there
exist elements $r_f\in \cam{E}p_{n_f}$ satisfying
$\sum_{i=1}^d\sum_{b\in\mathcal{B}^i}r^i_bb=\sum_{f\in\mathcal{F}'}r_ff$
(here $p_{n_f}\in P$ is such that $p_{n_f}f=f$). Since
$\mathcal{F}'\subseteq N_n$, $\mathcal{B}\cap M_n=\emptyset$ and
$\mathcal{F}'$ is $P$-$\mu$-independent, by the weak algorithm we
get an element $b'\in\mathcal{B}^{i'}\subseteq\mathcal{B}$ which is
$P$-$\mu$-dependent on
$(\mathcal{B}\setminus\{b'\})\cup\mathcal{F}'$. So, for all $i$, all
$b\in\mathcal{B}^i$ and all $f\in\mathcal{F}'$, there exist elements
$s^i_b\in \cam{E}p_i$, almost all zero, and elements $s_f\in
\cam{E}p_{n_f}$ such that
\[
\mu\Bigg(b'-\sum_{b\in\mathcal{B}\setminus\{b'\}}s_b^ib-\sum_{f\in\mathcal{F}'}s_ff\Bigg)<\mu(b')
\]
satisfying $\nu(s^i_b)+\mu(b)\leqslant\mu(b')$ and $\nu(s_f)+\mu
(f)\leqslant\mu(b')$. Moreover, we can assume that
$p_{i'}s^i_b=s^i_b$ and $p_{i'}s_f=s_f$. By the same argument used
in the proof of Proposition~\ref{prop:Ebases}, we see that
$\varphi(\mathcal{B}_{\mu(b')}^{i'})$ is linearly dependent modulo
$L'_{\mu(b')}$. This contradicts the fact that the image of
$\mathcal{B}_{\mu(b')}^{i'}$ is a $\cb$-basis of
$p_{i'}(L_{\mu(b')}/L'_{\mu(b')})$. Moreover, $M_n$ is
finite-dimensional over $\cb$ and $Q+\varphi(M_n)=L$ so $L/Q$ is
finite-dimensional over $\cb$.
\end{proof}

\begin{remark}\label{rem:fg}
Clearly, if $L$ in Theorem~\ref{tma:subproj} is finitely presented
then $Q$ is also finitely generated.
\end{remark}

\section{Flatness}

In this section, we prove that the Leavitt path algebra $\lev{E}$ is
flat as a {\em right} $\cam{\overline{E}}$-module. This will play an
important role in the sequel. We will denote by $\Si (E)$ the set of
vertices in $E$ which are sinks.

\begin{prop}\label{prop:flat}
$\lev{E}$ is flat as a right $\cam{\overline{E}}$-module.
\end{prop}
\begin{proof}
We write $R=\cam{\overline{E}}$ and $L=\lev{E}$. To prove that $L_R$
is flat, it suffices to show that $\mathrm{Tor}_1^R(L,M)=0$ for
every left $R$-module $M$. We will use the properties of quiver
algebras constructed in \cite[Section~2]{AraBrustenga07}. Recall
from there that the Leavitt path algebra is a quotient of
$S=(\cam{E})\left<\overline{E};\tau,\delta\right>$. More exactly,
let $X=E^0\setminus\Si (E) $ be the set of vertices which are not
sinks, then $L=S/I$, where $I$ is the ideal of $S$ generated by the
idempotent $q=\sum_{i\in X}p_i-\sum_{e\in E^1}e\overline{e}$ (see
\cite[Proposition~2.13]{AraBrustenga07}). From
\cite[Proposition~2.5]{AraBrustenga07} we know that elements in $S$
can be uniquely written as finite sums $\sum_{\alpha\in
E^*}r_\alpha\overline{\alpha}$, where $r_\alpha\in
\cam{E}p_{r(\alpha)}$. On the other hand, elements in $\cam{E}$ have
a unique expression as $\cb$-linear combinations of paths. We have
that
\begin{equation}\label{eq:SRpro}
S=\bigoplus_{\alpha\in
E^*}\cam{E}\overline{\alpha}=\bigoplus_{\substack{\alpha,\beta\in
E^*\\r(\alpha)=r(\beta)}}\cb
\beta\overline{\alpha}=\bigoplus_{\beta\in
E^*}\beta\left(\bigoplus_{\substack{\alpha\in
E^*\\r(\alpha)=r(\beta)}}\cb\overline{\alpha}\right)=\bigoplus_{\beta\in
E^*}\beta R;
\end{equation}
so, $S_R$ is projective.

Write $q_i=p_iqp_i$. Recall from the proof of (3) in \cite[Lemma~2.10]{AraBrustenga07}
that elements in $I$ can be uniquely written as finite sums
\[
\sum_{i\in X}\sum_{\{\alpha\in E^*\mid r(\alpha)=i\}}r_\alpha q_i
\overline{\alpha},
\]
where $r_\alpha\in \cam{E}p_{r(\alpha)}$. Thus,
% \[
% I=\bigoplus_{i\in X}\bigoplus_{\{\alpha\in E^*\mid s(\alpha)=i\}}\overline{\alpha}q_i\cam{E},
% \]
proceeding in the same way as in \eqref{eq:SRpro} we get that
\[
I=\bigoplus_{i\in X}\bigoplus_{\{\gamma\in E^*\mid
r(\gamma)=i\}}\gamma q_iR
\]
is projective as right $R$-module.

Now, the exact sequence of right $R$-modules
\[
0\longrightarrow I\longrightarrow S\longrightarrow L\longrightarrow 0,
\]
gives a projective resolution for $L$. Let $M$ be a left $R$-module.
We want to see that the induced homomorphism
\[
\varphi\colon \bigoplus_{i\in X}\bigoplus_{\{\gamma\in E^* \mid
r(\gamma)=i\}}\gamma q_i R \otimes_R M \cong I\otimes_R
M\longrightarrow S\otimes_R M\cong \bigoplus_{\gamma\in E^*}\gamma
R\otimes_R M
\]
is a monomorphism. We observe that
\[
\varphi\left(\sum_{i\in X}\sum_{\{\gamma\in E^*\mid
r(\gamma)=i\}}\gamma q_i \otimes m_{\gamma}\right)= \sum_{i\in
X}\sum_{\{\gamma\in E^*\mid r(\gamma)=i\}}\left(\gamma \otimes
m_\gamma -\!\sum_{e\in s^{-1}(i)}\gamma e \otimes
\overline{e}m_\gamma\right),
\]
and pick a non-zero element
\[
x=\sum_{i\in X}\sum_{\{\gamma\in E^*\mid r(\gamma)=i\}}\gamma q_i
\otimes m_{\gamma}\in \bigoplus_{i\in X}\bigoplus_{\{\gamma\in E^*
\mid r(\gamma)=i\}}\gamma q_iR\otimes M.
\]
Let $\gamma_0$ be a path of minimum length such that
$p_im_{\gamma_0}\neq 0$, where $i=r(\gamma_0)$. Since $\gamma
_0R\otimes_R M\cong p_iM$, we get $\gamma _0 \otimes
m_{\gamma_0}\neq 0$. Note also that the term $\gamma _0 \otimes
m_{\gamma_0}$ cannot be cancelled in $\varphi(x)$, because for each
of the non-zero terms $\gamma e \otimes \overline{e}m_\gamma$
appearing in that expression, the length of $\gamma e$ is strictly
larger than the length of $\gamma_0$, and the sum
$\bigoplus_{\gamma\in E^*}\gamma R\otimes_R M $ is a direct sum. It
follows that $\varphi$ is injective and so
$\mathrm{Tor}_1^R(L,M)=0$, as desired.
\end{proof}

As a consequence, we can regard Leavitt path algebras as perfect
left localizations (see \cite[Chapter~XI]{Stenstrom75}) of path
algebras:
\begin{cor}
The Leavitt path algebra $\lev{E}$ is a flat epimorphic left ring of
quotients of $\cam{\overline{E}}$.
\end{cor}
\begin{proof}
For $i\in E^0\setminus\Si (E) $, we write
$s^{-1}(i)=\{e^i_1,\dotsc,e^i_{n_i}\}$ and consider the left
$\cam{\overline{E}}$-module homomorphisms
\begin{align*}
\nu_i\colon\bigoplus_{j=1}^{n_i} \cam{\overline{E}} p_{s(e^i_j)} &\longrightarrow \cam{\overline{E}}p_i\\
(r_1,\dotsc,r_{n_i})&\longmapsto\sum_{j=1}^{n_i}r_j
\overline{e}^i_j.
\end{align*}
We write $\Sigma_2=\{\nu_i\mid i\in E^0\setminus\Si (E) \}$ (see the
Introduction). It is easy to see that the inclusion
$\cam{\overline{E}}\hookrightarrow \lev{E}$ is a universal
$\Sigma_2$-inverting homomorphism; so, it is a ring epimorphism (see
\cite[Chapter~4]{Schofield85}) and by Proposition~\ref{prop:flat} we
get that $\lev{E}$ is flat as a right $\cam{\overline{E}}$-module,
as desired.
\end{proof}

\begin{remark} \label{rem:quotients}

{\rm (1)} It is easy to see that the maximal flat epimorphic left
ring of quotients of $\cam{\overline{E}}$ is given by the regular
algebra of $E$, i.e. the algebra $\reg{E}$ defined in
\cite{AraBrustenga07}, see also the Introduction.

{\rm (2)} The fact that $\lev{E}$ is a left quotient ring of
$\cam{\ol{E}}$ (equivalently, a right quotient ring of $\cam{E}$)
has been already observed in \cite[Proposition 2.2]{Siles08}.

\end{remark}

\section{Finitely presented modules over the Leavitt path algebra}\label{sec:fpmL}
Recall that for every left semihereditary ring $S$, the category of
finitely presented left $S$-modules $\fp{S}$ is an abelian category.
(Here, we are looking at $\fp{S}$ as a full subcategory of the
category $\Mod{S}$ of all left $S$-modules. The fact that $S$ is
left semihereditary  implies that the kernel, the image and the
cokernel of every map between finitely presented modules are also
finitely presented).

We write $R=\cam{\overline{E}}$ for a finite quiver $E$, and let
$\mathcal{T}$ be the full subcategory of $\Mod{R}$ consisting of all
the left $R$-modules of finite dimension over $\cb$. This category
is obviously an abelian category, and we will show below that it is
the category of objects with finite length in the category $\fp{R}$.

\begin{prop}\label{prop:fdfl}
The category $\mathcal{T}$ of finite-dimensional left $R$-modules
coincides with the category $\fpl{R}$ of modules with finite length
in $\fp{R}$.
\end{prop}
\begin{proof}
First of all, note that every finite-dimensional left $R$-module is
finitely presented by Proposition~\ref{prop:LewinSchreier}. Clearly
all the objects in $\mathcal{T}$ are objects of finite length in
$\fp{R}$. It remains to see that a simple object in $\fp{R}$ must be
finite-dimensional. Let $M$ be a simple object in $\fp{R}$. By
Theorem~\ref{tma:subproj} (and Remark~\ref{rem:fg}), there is a
finitely generated projective $R$-module $Q$ such that $Q\leqslant
M$ and $M/Q$ is finite-dimensional. Since $M$ is simple in $\fp{R}$,
we must have $Q=0$; thus $M$ is finite-dimensional.
\end{proof}

We write $R=\cam{\overline{E}}$ for some finite quiver $E$, and
$A_E$ for the adjacency matrix of the quiver $E$.
\begin{prop}
Let $\mathcal{T}$ be the category of finite-dimensional left
$R$-modules. Then the following properties hold:
\begin{enumerate}
 \item $K_0(\mathcal{T})$ is a free abelian group over the set of isomorphism classes of simple,
 finite-dimensional left $R$-modules.\label{it:fag}
 \item The canonical map $\iota\colon K_0(R)\to K_0(\fp{R})$ is an isomorphism, so that $K_0(\fp{R})$
 is a free abelian group freely generated by $[Rp_1],\dotsc,[Rp_d]$.\label{it:iso}
 \item The map $K_0(\mathcal{T})\to K_0(\fp{R})$ sends $K_0(\mathcal{T})$ onto
 the subgroup of $K_0(\fp{R})$ generated by the columns of the matrix $\id-A_E^t$.\label{it:img}
\end{enumerate}
\end{prop}
\begin{proof}
\eqref{it:fag} Since the category $\mathcal{T}$ coincides with $\fpl{R}$ by
Proposition~\ref{prop:fdfl}, the result follows from the {\em Devissage} Theorem \cite[Theorem~3.1.8]{Rosenberg94}.

\eqref{it:iso} Since $R$ is a left hereditary ring, this is a
consequence of the Resolution Theorem
\cite[Theorem~3.1.13]{Rosenberg94}.

\eqref{it:img} We will denote by $[P]$ the class of a projective $R$-module $P$ in
$K_0(R)$ and by $\left<M\right>$ the class of a finitely presented $R$-module $M$
in $K_0(\fp{R})$. Moreover, we will identify $K_0(\cb^d)$ with $K_0(R)$ using the
isomorphism induced by the inclusion $\cb^d\hookrightarrow R$.

Now, let $M$ be a finite-dimensional $R$-module, by Proposition~\ref{prop:LewinSchreier} it admits a resolution
\[
0\longrightarrow P\longrightarrow Q\longrightarrow M\longrightarrow 0,
\]
where $P$ and $Q$ are finitely generated projective left $R$-modules. %so $\left<M\right>=\left<Q\right>-\left<P\right>$ in
%$K_0(\fp{R})$.
By the identification above and Proposition~\ref{prop:LewinSchreier} we get the equation
\[
\chi_{R}(M)=(\id-A_E^t)\chi_{\cb^d}(M)
\]
in $K_0(R)$. Moreover, since $\chi_R(M)=[Q]-[P]$ we get
\[
\left<M\right>=\left<Q\right>-\left<P\right>=\iota(\chi_R(M))=(\id-A_E^t)\iota(\chi_{\cb^d}(M))
\]
in $K_0(\fp{R})$. Therefore, the image of $K_0(\mathcal{T})$ is
contained in the subgroup generated by the columns of $(\id-A_E^t)$.

To see the reverse inclusion, remember that if $i\in E^0$ is not a
source then we have defined the left $R$-module homomorphisms
$\nu_i$. If $i\in E^0$ is a source we define $\nu_i$ as the zero
homomorphism $0\to Rp_i$. Now, the class
$\left<\coker(\nu_i)\right>$ in $K_0(\fp{R})$ coincides with the
$i$-th column of $(\id-A_E^t)$.
\end{proof}

Let $\mathcal{M}_\infty$ be the full subcategory of
$\Mod{\cam{\overline{E}}}$ with objects the modules $M$ such that
$\lev{E}\otimes_{\cam{\overline{E}}} M=0$. Moreover, we will write
$\mathcal{M}$ for the full subcategory of $\mathcal{M}_\infty$ given
by its finitely presented modules.

Recall that a \emph{Serre subcategory} of an abelian category $\mathcal{A}$ is
an abelian subcategory $\mathcal{B}$ which is closed under subobjects, quotients and extensions.
It is easy to see that the kernel of an exact functor between abelian categories is a
Serre subcategory ({\em cf.} \cite[Exercise~6.3.5]{BerrickKeating00}), hence
the category $\mathcal{M}_\infty$ is a Serre subcategory of $\Mod{\cam{\overline{E}}}$.
\begin{lem}\label{lem:catM}
Objects in the category $\mathcal{M}$ are finitely presented $\cam{\overline{E}}$-modules of finite length.
In fact, $\mathcal{M}$ is a Serre subcategory of $\fpl{\cam{\overline{E}}}$. Moreover,
the induced morphism $K_0(\mathcal{M})\to K_0(\fpl{\cam{\overline{E}}})$ is a monomorphism
and its image is the subgroup generated by the classes of simple modules in $\mathcal{M}$.
\end{lem}
\begin{proof}
Let $M$ be a module in $\mathcal{M}$. By Theorem~\ref{tma:subproj} $M$ has a (finitely generated)
projective submodule $P$ of finite codimension, so we have an exact sequence
\[
0\longrightarrow P\longrightarrow M\longrightarrow M/P\longrightarrow 0.
\]
Since ${\lev{E}}_{\cam{\overline{E}}}$ is flat
(Proposition~\ref{prop:flat}) we get that
$\lev{E}\otimes_{\cam{\overline{E}}}P=0$; hence $P=0$ and $M$ has
finite $\cb$-dimension. In particular $M$ has finite length.

We have an exact functor
$F\colon\fp{\cam{\overline{E}}}\to\fp{\lev{E}}$ given by
$F(M)=\lev{E}\otimes_{\cam{\overline{E}}}M$. It follows easily that
the kernel of this functor is precisely $\mathcal{M}$, thus
$\mathcal{M}$ is a Serre subcategory of both
$\fp{\cam{\overline{E}}}$ and $\fpl{\cam{\overline{E}}}$. Now, by
the {\em Devissage} Theorem (\cite[Theorem~3.1.8]{Rosenberg94}) we
are done.
\end{proof}

We shall need a result from \cite{Neeman07}. We have the following definition:
\begin{defi}[{\cite[Definition~0.4]{Neeman07}}]\label{def:catNee}
Let $R$ be a ring and let $\Sigma$ be a set of homomorphisms of finitely generated projective $R$-modules.
Assume all the maps in $\Sigma$ are monomorphisms. We define an exact category $\mathcal{E}$.
It is a full subcategory of all $R$-modules. All objects in $\mathcal{E}$ are
finitely presented $R$-modules, of projective dimension $\leqslant 1$.
The category $\mathcal{E}$ is completely determined by
\begin{enumerate}
 \item For every $s\colon P\to Q$ in $\Sigma$, the cokernel $M=Q/P$ lies in $\mathcal{E}$.\label{it:first}
 \item In any short exact sequence of finitely presented $R$-modules of projective dimension $\leqslant 1$
 \[
  0\longrightarrow M'\longrightarrow M\longrightarrow M''\longrightarrow 0,
 \]
 if two of the objects $M'$, $M$ and $M''$ lie in $\mathcal{E}$ then so does the third.
 \item $\mathcal{E}$ contains all direct summands of its objects.\label{it:last}
 \item $\mathcal{E}$ is minimal, subject to \eqref{it:first}--\eqref{it:last}.
\end{enumerate}
\end{defi}

There is an alternative characterization for this torsion category:
\begin{prop}[{\cite[Proposition~0.7]{Neeman07}}]\label{prop:Neeman}
An $R$-module $M$ belongs to $\mathcal{E}$ if and only if
\begin{enumerate}
 \item $M$ is finitely presented, and of projective dimension $\leqslant 1$.
 \item $R{\Sigma}^{-1}\otimes_R M=0=\mathrm{Tor}_1^R(R{\Sigma}^{-1},M)$.
\end{enumerate}
\end{prop}

Following \cite{Neeman07}, we shall refer to
$\mathcal{E}=\mathcal{E}(R,\Sigma)$ as the {\em category of
$(R,\Sigma)$-torsion modules}. An object of $\mathcal{E}$ will be a
$(R,\Sigma)$-torsion module. Using these results, we can
characterize the category $\mathcal{M}$:
\begin{theor}\label{thm:catM}
The category $\mathcal{M}$ is the full subcategory of
$\fpl{\cam{\overline{E}}}$ whose objects are the modules having all
composition factors in $\{\coker{\nu_i}\mid i\in E^0\setminus\Si (E)
\}$.
\end{theor}

\begin{proof}
Let $M$ be a module in $\mathcal{M}$. By definition, the module $M$
is a finitely presented $\cam{\overline{E}}$-module such that
$\lev{E}\otimes_{\cam{\overline{E}}}M=0$. Moreover, since
$\lev{E}_{\cam{\overline{E}}}$ is flat (Proposition~\ref{prop:flat})
and $\cam{\overline{E}}$ is a hereditary ring the remaining
conditions in Proposition~\ref{prop:Neeman} are fulfilled. Hence we
get $\mathcal M=\mathcal{E} (\cam{\ol{E}},\Sigma_2)$ from
Proposition \ref{prop:Neeman}.

Let $\mathcal M '$ be the category described in the statement. It is
clear that $\mathcal M'$ verifies (1)--(4) of Definition
\ref{def:catNee}. Thus we get
$$\mathcal M '=\mathcal{E} (\cam{\ol{E}},\Sigma_2)= \mathcal M \, ,$$
as desired.
\end{proof}

In order to obtain a description of the finitely presented $\lev{E}$-modules of
finite length we will need the following lemmas (\emph{cf.}~\cite[Lemma~6.1]{Sheiham06}):

\begin{lem}\label{lem:simple}
Let $N$ be a finite-dimensional simple $\cam{\overline{E}}$-module.
We have the following dichotomy:
\begin{enumerate}
 \item There exist $i\in E^0\setminus\Si (E) $ such that $N\cong\coker\nu_i$.
 In this case $\lev{E}\otimes_{\cam{\overline{E}}}N=0$.\label{item1simple}
 \item For every $i\in E^0\setminus\Si (E)$ we have $N\not\cong\coker\nu_i$.
 In this situation $\lev{E}\otimes_{\cam{\overline{E}}} N $ is simple.\label{item2simple}\end{enumerate}
\end{lem}

\begin{proof}
\eqref{item1simple} If $N\cong\coker\nu_i$ for some $i$ then
$\lev{E}\otimes_{\cam{\overline{E}}}N=0$ because
$\coker\nu_i\in\mathcal{M}$.

\eqref{item2simple} Let $N$ be a finite-dimensional simple left
$\cam{\overline{E}}$-module such that, for every $i\in E^0\setminus
\Si (E)$, we have $N\not\cong\coker\nu_i$. Theorem~\ref{thm:catM}
implies that $N\notin \mathcal{M}$, so that
$\lev{E}\otimes_{\cam{\overline{E}}}N\ne 0$.

Let $n=\sum_{\gamma\in E^*}\gamma \otimes n_{\gamma} $ be a nonzero
element in $\lev{E}\otimes_{\cam{\overline{E}}}N$, where
$n_{\gamma}\in N$. We may consider the following decomposition of
the unit
\[
1=\sum_{i\in E^0}p_i=\sum_{i\not\in\Si (E)}\sum_{e\in
s^{-1}(i)}e\overline{e}+\sum_{i\in\Si (E)}p_i.
\]
If $n':=p_in\neq 0$ for some sink $i$ then $n'\in p_i\otimes
p_iN\subseteq 1\otimes N$. Otherwise, we see that there is some
$e\in E^1$ such that $\overline{e}n\neq 0$, and we see inductively
that we can find $\gamma\in E^*$ such that
$n':=\overline{\gamma}n\neq 0$ and $n'\in 1\otimes N$. In both
cases, the simplicity of $N$ gives us $\cam{\overline{E}}n'=1\otimes
N$, showing the simplicity of
$\lev{E}\otimes_{\cam{\overline{E}}}N$.
\end{proof}

\begin{lem}\label{lem:fl}
Let $i$ be a vertex. The following are equivalent:
\begin{enumerate}
 \item $\cam{\overline{E}}p_i$ has finite $\cb$-dimension.\label{item1:fl}
 \item $\lev{E}p_i$ is a finite direct sum of simple submodules.\label{item2:fl}
 \item $\lev{E}p_i$ has finite length.\label{item3:fl}
 \item The subgraph $s^{-1}_{E^*}(i)$ is acyclic.\label{item4:fl}
\end{enumerate}
\end{lem}
\begin{proof}
\eqref{item1:fl}$\Rightarrow$\eqref{item2:fl}. Let
$M\subseteq\cam{\overline{E}}$ be the set of all paths in
$\overline{E}$ with range $i$ and starting at a source of
$\overline{E}$. Since $\cam{\overline{E}}p_i$ has finite
$\cb$-dimension, the set $M$ is finite. We remark that every path in
$\overline{E}$ with range $i$ can be extended to a path in $M$. Now,
using the relations $p_j=\sum_{e\in s^{-1}(j)}e\overline{e}$
iteratively and the previous remark we get that
$p_i=\sum_{\overline{\gamma}\in M}\gamma \overline{\gamma}$, hence
$\lev{E}p_i=\sum_{\overline{\gamma}\in M}\lev{E}\gamma
\overline{\gamma}$. Moreover, this is a direct sum because elements
in the set $\{\gamma \overline{\gamma} \mid\overline{\gamma}\in M\}$
are orthogonal idempotents. On the other hand,
\[
\lev{E}\gamma
\overline{\gamma}\cong\lev{E}\overline{\gamma}\gamma=\lev{E}p_{r(\gamma)}\cong
\lev{E}\otimes_{\cam{\overline{E}}}\cam{\overline{E}}p_{r(\gamma)}.
\]
Since, $r(\gamma)=s(\overline{\gamma})$ is a source in
$\overline{E}$ the module $\cam{\overline{E}}p_{r(\gamma)}$ is
simple and we are done by Lemma~\ref{lem:simple}.

\eqref{item2:fl}$\Rightarrow$\eqref{item3:fl} is obvious.

\eqref{item3:fl}$\Rightarrow$\eqref{item4:fl}. Suppose that the
subgraph $s^{-1}_{E^*}(i)$ is not acyclic. In particular, there are
paths $\alpha,\gamma\in E^*$ such that $\alpha$ is a cycle based at
some vertex $k$, $r(\gamma)=k$ and $s(\gamma)=i$. We write
$x=p_i+\gamma \alpha\overline{\gamma}$. If $n>m\geqslant 1$ are
natural numbers then $\lev{E}x^n\subset \lev{E}x^m$. Indeed, suppose
$y\in\lev{E}$ is such that $yx^n=x^m$. Since $p_i\lev{E}p_i\subseteq
p_i\reg{E}p_i$, operating in the latter ring we get that
$y=x^{m-n}$, but $m-n<0$ and hence $y\not\in p_i\lev{E}p_i$.
Therefore, we have constructed an infinite chain of submodules with
proper inclusions:
\[
\lev{E}x\supset \lev{E}x^2\supset \cdots\supset
\lev{E}x^n\supset\cdots
\]

\eqref{item4:fl}$\Rightarrow$\eqref{item1:fl} is clear.
\end{proof}
Our next result gives a description of the structure of the finitely
presented $\lev{E}$-modules.

\begin{prop}\label{prop:modfp}
Let $E$ be a finite quiver and write $R=\cam{\overline{E}}$,
$L=\lev{E}$. Then the following holds:
\begin{enumerate}
 \item\label{item1:modfp} Let $N$ be a finite-dimensional left $R$-module with a composition series of length $k$:
  \[
   0< N_1< N_2<\dotsc< N_k=N.
  \]
 Assume that exactly $r$ composition factors are isomorphic to modules in
 the set $\{\coker\nu_i\mid i\in E^0\setminus\Si (E) \}$. Then $L\otimes_RN$ is a left
  $L$-module of finite length and its length is exactly $k-r$.
 \item\label{item2:modfp} Let $M$ be a finitely presented left $L$-module.
 Then there is a finitely generated projective $L$-module $P$ such that $P\leqslant M$
 and $M/P$ is a module of finite length.
 \item\label{item3:modfp} Every finitely presented left $L$-module $M$ of finite length
 is isomorphic to a module of the form $L\otimes_RN$, where $N$ is a finite-dimensional left
  $R$-module.
\end{enumerate}
\end{prop}

\begin{proof}
\eqref{item1:modfp} It follows easily from Lemma~\ref{lem:simple}
and the fact that $L$ is flat as a right $R$-module
(Proposition~\ref{prop:flat}).

\eqref{item2:modfp} Let $M$ be a finitely presented left $L$-module.
By \cite[Corollary~4.5]{Schofield85} there exists a finitely
presented left $R$-module $N$ such that $L\otimes_RN\cong M$. Now,
by Theorem~\ref{tma:subproj} (and Remark~\ref{rem:fg}), there is a
finitely generated projective $R$-module $Q$ such that $Q\leqslant
N$ and $N/Q$ is finite-dimensional. Since $L_R$ is flat, we have
that $M\cong L\otimes_RN$ contains the f.g. projective $L$-module
$P\cong L\otimes_RQ$. By \eqref{item1:modfp}, the $L$-module
$(L\otimes_RN)/(L\otimes_RQ)\cong L\otimes _R (N/Q)$ is of finite
length.

\eqref{item3:modfp} As above we know that $M\cong L\otimes_RN$ for
some finitely presented left $R$-module $N$ and we obtain (by
Theorem~\ref{tma:subproj}) a projective left $R$-module $Q$ such
that $N/Q$ is finite-dimensional. From the following exact sequence
\[
0\longrightarrow L\otimes_RQ\longrightarrow M \longrightarrow
L\otimes _R(N/Q)\longrightarrow 0
\]
we get that the projective left $L$-module $L\otimes_RQ$ has finite
length. Since $Q\cong\oplus_{i=1}^k Rp_{j_i}$ for some $j_i\in E^0$
and every $L\otimes _R Rp_{j_i} \cong Lp_{j_i}$ has finite length,
by Lemma~\ref{lem:fl} we get that every $Rp_{j_i} $ is
finite-dimensional. Thus, $Q$ is also finite-dimensional, and
therefore so is $N$.
\end{proof}

\section{The category of finitely presented modules as a quotient category}
\label{sect:quotient}

In this section we will prove that the categories $\Mod{\lev{E}}$,
$\fp{\lev{E}}$ and $\fpl{\lev{E}}$ are equivalent, respectively,
to the quotient categories $\Mod{\cam{\overline{E}}}/\mathcal{M}_\infty$,
$\fp{\cam{\overline{E}}}/\mathcal{M}$ and $\fpl{\cam{\overline{E}}}/\mathcal{M}$.
The following results generalize \cite[Section~5]{Ara04} to the quiver setting,
although quite often the ideas behind the proofs follow \cite{Sheiham06},
where the similar case of the free group algebra is considered.

We first recall some basics on categories. Given a Serre subcategory $\mathcal{B}$
of an abelian category $\mathcal{A}$, one can consider a quotient abelian category
$\mathcal{A}/\mathcal{B}$ and an exact functor $T\colon\mathcal{A}\to\mathcal{A}/\mathcal{B}$
with the following universal property: given an exact functor $S\colon\mathcal{A}\to \mathcal{C}$
from $\mathcal{A}$ to an abelian category $\mathcal{C}$ such that $S(B)\cong 0$ for every
object $B$ of $\mathcal{B}$, there is a unique exact functor
$S'\colon\mathcal{A}/\mathcal{B}\to \mathcal{C}$ such that $S=S'T$ (see \cite[Chapter~II]{WeibelXX}).
If the category $\mathcal{A}$ is well-powered (that is, every object in $\mathcal{A}$ has
a set of representative subobjects) then we can assure the existence of the quotient
category $\mathcal{A}/\mathcal{B}$ for any Serre subcategory $\mathcal{B}$
(see \cite[Theorem~I.2.1]{Swan68}). Since we only deal with module categories this condition is always fulfilled.

Recall that, given a category $\mathcal{C}$ and a collection
$\Sigma$ of morphisms in $\mathcal{C}$, the \emph{localization of
$\mathcal{C}$ with respect to $\Sigma$} is a category
$\mathcal{C}_\Sigma$, together with a functor
$L\colon\mathcal{C}\to\mathcal{C}_\Sigma$ such that
\begin{enumerate}
 \item For every $s\in\Sigma$, $L(s)$ is an isomorphism.
 \item If $F\colon\mathcal{C}\to\mathcal{D}$ is any functor sending $\Sigma$ to isomorphisms in
 $\mathcal{D}$, then $F$ factors uniquely through $L\colon\mathcal{C}\to\mathcal{C}_{\Sigma}$.
\end{enumerate}
It turns out that the quotient category $\mathcal{A}/\mathcal{B}$
can also be obtained by localization of $\mathcal{A}$ with respect
to the collection of all \emph{$\mathcal{B}$-isos}, that is, those
maps $f$ such that $\ker(f)$ and $\coker(f)$ are in $\mathcal{B}$
(for details see \cite[Appendix in Chapter~II]{WeibelXX}). Thus, we
can make use of both universal properties for the quotient category.
Moreover, maps in $\mathcal{A}/\mathcal{B}$ are given by equivalence
classes $[(f,g)]$ of diagrams in $\mathcal{A}$,
\[
A_1\xleftarrow{\;f\;} A\xrightarrow{\;g\;} A_2
\]
where $f$ is a $\mathcal{B}$-iso.

Let us write $B=\lev{E}\otimes_{\cam{\overline{E}}} -
\colon\Mod{\cam{\overline{E}}}\to \Mod{\lev{E}}$ for the functor
given by extension of scalars and
$U\colon\Mod{\lev{E}}\to\Mod{\cam{\overline{E}}}$ for the functor
given by restriction of scalars. We remark that $B$ and $U$ are
adjoint functors (see \cite[Proposition~3.3.15]{BerrickKeating00}).
We know that $B$ restricts to a functor between the categories of
finitely presented modules and, by
Proposition~\ref{prop:modfp}\eqref{item1:modfp}, the same applies to
the subcategories of finite length modules. We will also denote
these restrictions by $B$.

Recall from Section~\ref{sec:fpmL} that $\mathcal{M}_\infty$ is a
Serre subcategory of $\Mod{\cam{\ol{E}}}$ and that $\mathcal{M}$ is
a Serre subcategory of $\fp{\cam{\ol{E}}}$ and of
$\fpl{\cam{\ol{E}}}$ (see Lemma~\ref{lem:catM}). Therefore, it makes
sense to consider the quotient categories
$\Mod{\cam{\overline{E}}}/\mathcal{M}_\infty$,
$\fp{\cam{\overline{E}}}/\mathcal{M}$ and
$\fpl{\cam{\overline{E}}}/\mathcal{M}$.

\begin{prop}\label{prop:pronat}
Let $M\in\Mod{\cam{\overline{E}}}$ and $N\in\Mod{\lev{E}}$. Then the following properties hold:
\begin{enumerate}
 \item\label{item1:pronat} There is a natural isomorphism $\eta_N\colon BU(N)\to N$.
 \item\label{item2:pronat} There is a natural transformation $\theta_M\colon M\to UB(M)$.
 \item\label{item3:pronat} The composites
\[
U(N)\xrightarrow{\theta_{U(N)}}{}UBU(N)\xrightarrow{U(\eta_N)}{}U(N)
\]
\[
B(M)\xrightarrow{B(\theta_M)}{}BUB(M)\xrightarrow{\eta_{B(M)}}{}B(M)
\]
are identity morphisms.
\end{enumerate}
\end{prop}
\begin{proof}
\eqref{item1:pronat} Recall that the inclusion
$\cam{\overline{E}}\hookrightarrow \lev{E}$ is a universal
localization; thus it is a ring epimorphism and, by
\cite[Proposition~XI.1.2]{Stenstrom75}, the natural transformation
$\eta_N\colon BU(N)\to N$ defined by $\eta_N(s\otimes n)=sn$ is a
natural isomorphism.

\eqref{item2:pronat} It is clear that the homomorphism
$\theta_M\colon M\to UB(M)$ defined by $\theta_M(m)=1\otimes m$ is
natural.

\eqref{item3:pronat} It is obvious from the previous definitions.
\end{proof}

We deduce in the next proposition that $B$ satisfies the same universal property as the localization functor,
but only up to natural isomorphism. Let $\Xi$ be the collection of all $\mathcal{M}_\infty$-isos in
$\Mod{\cam{\overline{E}}}$.
\begin{prop}\label{prop:natpro}
If $S\colon\Mod{\cam{\overline{E}}}\to\mathcal{B}$ is a functor which sends every
morphism in $\Xi$ to an isomorphism then there is a functor $S'\colon\Mod{\lev{E}}\to\mathcal{B}$
such that $S'B$ is naturally isomorphic to $S$. Moreover, the functor $S'$ is unique up to natural isomorphism.
\end{prop}
\begin{proof}
We prove uniqueness first. If there is a natural isomorphism
$S\simeq S'B$ then $SU\simeq S'BU\simeq S'$ by Proposition~\ref{prop:pronat}\eqref{item1:pronat}.

To prove existence we must show that if $S'=SU$ then $S'B\simeq S$. Indeed,
by Proposition~\ref{prop:pronat}\eqref{item3:pronat} $B(\theta_M)\colon B(M)\to BUB(M)$
is an isomorphism for each $M\in\Mod{\cam{\overline{E}}}$. Since $B$ is an exact functor
(Proposition~\ref{prop:flat}) we have $\theta_M\in\Xi$. Thus, $S(\theta)\colon S\to SUB=S'B$ is a natural isomorphism.
\end{proof}

Let us consider the localization functor:
\[
T\colon\Mod{\cam{\overline{E}}}\to\Mod{\cam{\overline{E}}}/\mathcal{M}_\infty.
\]
By the universal property of $T$ there exists a unique functor
\[
\overline{B}\colon\Mod{\cam{\overline{E}}}/\mathcal{M}_\infty\longrightarrow\Mod{\lev{E}}
\]
such that $B=\overline{B}T$. We will denote by $\fpl{\cam{\overline{E}}}/\mathcal{M}_\infty$ and $\fp{\cam{\overline{E}}}/\mathcal{M}_\infty$ the full subcategories of $\Mod{\cam{\overline{E}}}/\mathcal{M}_\infty$ given, respectively, by the finitely presented modules of finite length and by the finitely presented modules. Beware that $\mathcal{M}_\infty$ is not contained in the categories of finitely presented modules so, despite of the notation, these are not quotient categories.

We have the following commutative diagram:
\[
\xymatrix{
\fpl{\cam{\overline{E}}}\ar[r]^{T_{\rm fl}}\ar[d]&\frac{\fpl{\cam{\overline{E}}}}{\mathcal{M_\infty}}\ar[d]\ar[r]^{\overline{B}_{\rm fl}}&\fpl{\lev{E}}\ar[d]\\
\fp{\cam{\overline{E}}}\ar[r]^{T_{\rm fp}}\ar[d]&\frac{\fp{\cam{\overline{E}}}}{\mathcal{M_\infty}}\ar[d]\ar[r]^{\overline{B}_{\rm fp}}&\fp{\lev{E}}\ar[d]\\
\Mod{\cam{\overline{E}}}\ar[r]^{T}\ar@/_1.5pc/[rr]_B&\frac{\Mod{\cam{\overline{E}}}}{\mathcal{M}_\infty}\ar[r]^{\overline{B}}&\Mod{\lev{E}} }
\]
where the vertical arrows are inclusions of full subcategories and the horizontal ones in the first and second rows are given by restriction.

\begin{theor}\label{tma:eqscat}
The functors $\overline{B}$, $\overline{B}_{\rm fp}$ and $\overline{B}_{\rm fl}$ are category equivalences.
\end{theor}
\begin{proof}
Recall that two categories are equivalent if and only if there is a full, faithful and dense functor between them (see \cite[Proposition~1.3.14]{BerrickKeating00}). By Proposition~\ref{prop:natpro}, the functor $B$ satisfies the same natural property than $T$ up to natural isomorphism, hence $\overline{B}$ is a category equivalence. Since $\overline{B}_{\rm fp}$ and $\overline{B}_{\rm fl}$ are given by restriction of $\overline{B}$, these are full and faithful functors. Moreover, the functor $\overline{B}_{\rm fp}$ is dense by \cite[Corollary~4.5]{Schofield85} and $\overline{B}_{\rm fl}$ is dense as a consequence of Proposition~\ref{prop:modfp}\eqref{item3:modfp}.
\end{proof}

\begin{prop}\label{prop:quoeqv}
The following holds:
\begin{enumerate}
\item The category $\fpl{\cam{\overline{E}}}/\mathcal{M_\infty}$ is equivalent to the quotient category $\fpl{\cam{\overline{E}}}/\mathcal{M}$.
\item The category  $\fp{\cam{\overline{E}}}/\mathcal{M_\infty}$ is equivalent to the quotient category  $\fp{\cam{\overline{E}}}/\mathcal{M}$.
\end{enumerate}
\end{prop}
\begin{proof}
Let us consider the localization functor in each case:
\begin{align*}
S_{\rm fl}\colon\fpl{\cam{\overline{E}}}&\longrightarrow\fpl{\cam{\overline{E}}}/\mathcal{M}\\
S_{\rm fp}\colon\fp{\cam{\overline{E}}}&\longrightarrow\fp{\cam{\overline{E}}}/\mathcal{M}.
\end{align*}
By the universal property there exist two unique functors
\begin{align*}
\overline{T}_{\rm fl}\colon\fpl{\cam{\overline{E}}}/\mathcal{M}&\longrightarrow\fpl{\cam{\overline{E}}}/\mathcal{M_\infty}\\
\overline{T}_{\rm fp}\colon\fp{\cam{\overline{E}}}/\mathcal{M}&\longrightarrow\fp{\cam{\overline{E}}}/\mathcal{M_\infty}
\end{align*}
satisfying that $T_{\rm fl}=\overline{T}_{\rm fl}S_{\rm fl}$ and $T_{\rm fp}=\overline{T}_{\rm fp}S_{\rm fp}$. We will show that $\overline{T}_{\rm fp}$ is a full, faithful and dense functor, hence a category equivalence. %We will only comment the differences for the similar case of $\overline{T}_{\rm fl}$.

Since the categories $\fp{\cam{\overline{E}}}/\mathcal{M}$ and $\fp{\cam{\overline{E}}}/\mathcal{M_\infty}$ have the same objects and $\overline{T}_{\rm fp}$ acts as the identity on them it is a dense functor in a trivial way.

Let us write $F=\overline{B}_{\rm fp}\overline{T}_{\rm fp}$. The
maps in $\fp{\cam{\overline{E}}}/\mathcal{M}$ are equivalence
classes $[(f,g)]$ of diagrams in $\fp{\cam{\overline{E}}}$,
\[
 M_1\xleftarrow{\;f\;}M\xrightarrow{\;g\;}M_2
\]
where the kernel and the cokernel of $f$ are objects in
$\mathcal{M}$. For such a pair, we have $F([(f,g)])=(\id\otimes
g)(\id\otimes f)^{-1}$. Now assume that $(\id\otimes g)(\id\otimes
f)^{-1}=0$. Then $\id\otimes g =0$, so
$\im(g)\in\mathcal{M}_\infty$. Since $\fp{\cam{\ol{E}}}$ is an
abelian category and $\im(g)=\ker(\coker(g))$ this module is
finitely presented and hence in $\mathcal{M}$. Consequently
$[(f,g)]=[(f,0)]=0$ and $F$ is a faithful functor. Therefore
$\overline{T}_{\rm fp}$ is faithful as well.

Now we will prove that $\overline{T}_{\rm fp}$ is  a full functor.
Let $M_1$ and $M_2$ be finitely presented right $\cam{\overline{E}}$-modules.
A map in $\fp{\cam{\overline{E}}}/\mathcal{M}_{\infty}$ is given by an equivalence
class $[(f,g)]$ of diagrams in $\Mod{\cam{\overline{E}}}$,
\[
 M_1\xleftarrow{\;f\;}M\xrightarrow{\;g\;}M_2
\]
where $M$ is a left $\cam{\overline{E}}$-module and the kernel and
the cokernel of $f$ are objects in $\mathcal{M}_{\infty}$. It is
enough to show that it is possible to pick a representative of
$[(f,g)]$ with $M$ finitely presented.

Let us write $N'=(\ker f)\cap(\ker g)$. From the following commutative diagram:
\[
\xymatrix{
& M\ar[ld]_f\ar[rd]^g\ar[d]^{\pi'} & \\
M_1 & M/N'\ar[r]^{\overline{g}}\ar[l]_{\overline{f}} & M_2
}
\]
we obtain that $[(f,g)]=[(\overline{f},\overline{g})]$. So we can assume that $f\oplus g\colon M\to M_1\oplus M_2$ is a monomorphism.

We will show that for such an $M$ we have $M\in\fp{\cam{\overline{E}}}$. By Theorem~\ref{tma:subproj} (and Remark~\ref{rem:fg}) there exist finitely generated and projective submodules $P_1\subseteq M_1$, $P_2\subseteq M_2$ such that $M_1/P_1$ and $M_2/P_2$ have finite dimension.
Let us write $\pi_1\colon M_1\to M_1/P_1$ and $\pi_2\colon M_2\to M_2/P_2$ for the natural projections and consider the module
\[
N=(\ker\pi_1f)\cap(\ker\pi_2g).
\]

We have the following commutative diagram with exact rows:
\begin{equation}\label{eq:resnopro}
\begin{CD}
%@. @. 0 @. 0 \\
%@. @. @VVV @VVV\\
0 @>>> N @>>> M @>>> M/N @>>> 0\\
@. @VVf'\oplus g'V @VVf\oplus gV @VVf''\oplus g''V \\
0 @>>> P_1\oplus P_2 @>>> M_1\oplus M_2 @>>> M_1/P_1\oplus M_2/P_2@>>> 0,
\end{CD}
\end{equation}
where $f'\oplus g'$ is induced by the universal property of the kernel and $f''\oplus g''$ is induced by the universal property of the cokernel. Observe that the vertical arrows are monomorphisms. Therefore the module $N$ is projective and the module $M/N$ has finite dimension (and by Proposition~\ref{prop:fdfl} is finitely presented).

Consider a resolution of $M/N$ by finitely generated projective $\cam{\overline{E}}$-modules:
\[
 0\longrightarrow Q\longrightarrow P\longrightarrow M/N\longrightarrow 0.
\]
Applying Schanuel Lemma (\cite[(5.1)]{Lam99}) to the previous resolution and to the first row in \eqref{eq:resnopro} we get the following projective resolution of $M$:
\[
0\longrightarrow Q\longrightarrow N\oplus P \longrightarrow M\longrightarrow 0.
\]

We just need to check that $N\oplus P$ is finitely generated. Recall
that in a semihereditary ring every projective module is isomorphic
to a direct sum of finitely generated ideals (see
\cite[Theorem]{Albrecht61}). Thus, we may consider the following
decomposition into direct summands $N\oplus P=Q_1\oplus Q_2$, where
$Q\subseteq Q_1$ and $Q_1$ is a finitely generated projective
module. Now $M\cong (Q_1/Q)\oplus Q_2$ decomposes as a direct sum of
a projective module and a finitely presented module. We obtain
\[
\lev{E}\otimes_{\cam{\overline{E}}} M_1\cong
\lev{E}\otimes_{\cam{\overline{E}}}M \cong
 \left(\lev{E}\otimes_{\cam{\overline{E}}}(Q_1/Q)\right)\oplus\left(\lev{E}\otimes_{\cam{\overline{E}}}Q_2
\right).
\]
Since the module $\lev{E}\otimes_{\cam{\overline{E}}}M_1$ is
finitely presented, the module
$\lev{E}\otimes_{\cam{\overline{E}}}Q_2$ is finitely presented as
well. Now, since $Q_2$ is projective, we get that $Q_2$ is finitely
generated and $M$ is finitely presented. Moreover,
$\ker(f),\coker(f)\in\mathcal{M}$ and we have seen that the functor
$\overline{T}_{\mathrm{fp}}$ is full.

The proof for $\overline{T}_{\mathrm{fl}}$ is similar, but simpler because $\fpl{\cam{\overline{E}}}$ is closed under subobjects.
\end{proof}

As a consequence of Theorem~\ref{tma:eqscat} and Proposition~\ref{prop:quoeqv} we obtain:
\begin{cor}
The following holds:
\begin{enumerate}
\item The categories $\fpl{\cam{\overline{E}}}/\mathcal{M}$ and $\fpl{\lev{E}}$ are equivalent.
\item The categories  $\fp{\cam{\overline{E}}}/\mathcal{M}$ and  $\fp{\lev{E}}$ are equivalent.
\end{enumerate}
\end{cor}

\section{Blanchfield modules over a quiver}\label{sec:Blanchfield}

Let $R$ be a ring and let $\Sigma$ be a family of injective
homomorphisms between finitely generated projective $R$-modules.
Recall that, by \cite[Proposition 2.2]{Neeman07}, all maps in
$\Sigma$ are injective in case the localization  map $R\to
R\Sigma^{-1}$ is injective.

The localization $R\to R\Sigma^{-1}$ is {\em stably flat} if
$\text{Tor}^R_i(R\Sigma^{-1},R\Sigma^{-1})=0$ for all $i\ge 2$.
Observe that if $R$ is left hereditary then every universal
localization $R\to R\Sigma^{-1}$ is stably flat. Moreover by a
result of Bergman and Dicks \cite[Theorem 5.3]{BergmanDicks78},
$R\Sigma^{-1}$ is also left hereditary.

\begin{theor}[Neeman, Ranicki \cite{NeemanRanicki01},\cite{NeemanRanicki04},\cite{Neeman07}]
\label{theor:Nee-Ran} Let $R\to R\Sigma^{-1}$ be a stably flat
universal localization such that all the morphisms in $\Sigma$ are
injective. Then there is an exact sequence in nonnegative $K$-theory
$$\cdots \to K_{i+1}(R)\to K_{i+1}(R\Sigma^{-1})\to K_i(\mathcal{E}(R,\Sigma))\to K_i(R)\to \cdots .$$
\end{theor}

Following terminology suggested by \cite{RanickiSheiham06}, we call
a left module $M$ over $\cam{E}$ a {\it Blanchfield module} in case
$\text{Tor}^{\cam{E}}_q(\cb^d, M)=0$  for all $q$, where we see
$\cb^d$ as a right $\cam{E}$-module through the augmentation
$\epsilon \colon \cam{E}\to \cb^d$. It is easy to check that $M$ is
a Blanchfield module if and only if the natural map
$$\bigoplus_{e\in E^1} p_{r(e)}M\longrightarrow M, \qquad (p_{r(e)}m_e)\mapsto \sum _{e\in E^1} em_e$$
is an isomorphism (see the proof of Proposition \ref{torsions} for
details). Note that this is equivalent to saying that $p_iM=0$ for
every $i\in \Si (E)$ and that all the maps $\bigoplus_{e\in
s^{-1}(i)} p_r(e)M\longrightarrow p_iM$, for $i\in E^0\setminus \Si
(E)$, are isomorphisms. It follows that the Blanchfield modules are
exactly the left $\lev{E}$-modules $M$ such that $p_iM=0$ for every
$i\in \Si (E)$.

We will denote the full subcategory of $\Mod{P(E)}$ consisting of
all the Blanchfield $P(E)$-modules by $\Bla _{\infty}(P(E))$, and
the category of finitely generated Blanchfield $P(E)$-modules by
$\Bla (P(E))$. Let $M$ be a f.g. Blanchfield $\cam{E}$-module. A
{\it lattice} in $M$ is a $\cam{\ol{E}} $-submodule $A\subset M$
such that $A$ is finite dimensional over $\cb$ and $M=\cam{E} A$.

For a ring $R$, denote by $\fnp{R}$ the full subcategory of finitely
presented $R$-modules of finite length without nonzero projective
submodules.

\begin{prop}
\label{typeLL}

 {\rm (1)} Let $M$ be a left $\lev{E}$-module. Then $M$ is a f.g. Blanchfield
 $\cam{E}$-module if and only if $M\in \fnp{L(E)}$.

 {\rm (2)} Let $M$ be a f.g. Blanchfield $\cam{E}$-module. Then $M$ contains a
 lattice. Moreover a $\cam{\ol{E}}$-submodule $A$ of $M$ is a lattice if and
only if $A$ is finite dimensional and the natural map
$\lev{E}\otimes _{{\cam{\ol{E}}}} A\to M$ is an isomorphism.
Furthermore, any lattice in $M$ does not contain nonzero projective
$\cam{\ol{E}}$-submodules.

{\rm (3)} Every f.g. Blanchfield $\cam{E}$-module contains a
smallest lattice.

\end{prop}

\begin{proof}
(1) If $M$ is a finitely presented $\lev{E}$-module of finite length
without nonzero projective submodules then by Proposition
\ref{prop:modfp}(3) there is a finite dimensional left
$\cam{\ol{E}}$-module $N$ such that $L(E)\otimes
_{\cam{\ol{E}}}N\cong M$. Then clearly $M$ is finitely generated as
a $\cam{E}$-module. If $i\in \Si (E)$ and $p_iM\ne 0$, then there is
a nonzero map $\lev{E}p_i\to M$ which is injective because
$\lev{E}p_i$ is simple, contradicting the fact that $M$ does not
contain nonzero projective submodules.

The converse follows from (2).

(2) Assume that $M$ is a left $\lev{E}$-module which is finitely
generated as $\cam{E}$-module. Let $a_1,\dots ,a_r$ generators of
$M$ as a left $\cam{E}$-module. Then, for $e\in E^1$,
$$\ol{e}a_i=\sum _k \gamma _{ji}^e a_j $$
where $\gamma_{ji}^e\in \cam{E}$. Let $r$ be an upper bound for the
lengths of the paths involved in the $\gamma_{ji}^e$'s. Let $A$ be
the $\cb$-space generated by $\lambda a_j$, where $|\lambda |\le r$.
Then $\ol{e}\lambda a_j\in A$, and clearly $A$ is a lattice for $M$.

If $A\subset M$ is a finite-dimensional $\cam{\ol{E}}$-submodule and
the natural map $\lev{E}\otimes _{\cam{\ol{E}}}A\rightarrow M$ is an
isomorphism, then $M=\cam{E} A$ and thus $A$ is a lattice in $M$.
Conversely assume that $A$ is a lattice in $M$. Since $\lev{E}$ is
flat as a right $\cam{\ol{E}}$-module, the map $\lev{E}\otimes
_{\cam{\ol{E}}}A\rightarrow \lev{E}\otimes _{\cam{\ol{E}}}M$ is
injective. Now the natural map $\lev{E}\otimes
_{\cam{\ol{E}}}M\rightarrow M$ is an isomorphism, because the
inclusion $\cam{\ol{E}}\rightarrow \lev{E}$ is a ring epimorphism.
It follows that the map $\lev{E}\otimes_{\cam{\ol{E}}} A\rightarrow
M$ is injective. Since $A$ is a lattice this map is clearly
surjective.

It follows that $M$ is a finitely presented $\lev{E}$-module of
finite length. If $p_iM=0 $ for every $i\in \Si (E)$ then $M$ does
not have nonzero projective submodules by Lemma \ref{lem:fl}.
Observe that this implies that any lattice $A$ of $M$ does not
contain nonzero projective $\cam{\ol{E}}$-submodules.

(3) This follows as in \cite[Proposition 4.1(3)]{Ara04}, by showing
that the intersection of two lattices is a lattice.
\end{proof}

Let $\Sigma $ be the set of square matrices over $\cam{E}$ that are
sent to invertible matrices by the augmentation homomorphism
$\epsilon\colon \cam{E}\to \cb^d$. We have $\rat{E}\cong
\cam{E}\Sigma ^{-1}$, see diagram (\ref{maindiagram}) and the
comments below it. We are now ready to determine the categories of
$(\cam{E},\Sigma)$-torsion and $(\lev{E},\Sigma)$-torsion.

\begin{prop}
\label{torsions} With the above notation, we have
$$\mathcal{E}(\cam{E},\Sigma)=\Bla (\cam{E})=\mathcal{E}(\lev{E},\Sigma).$$
Moreover $\Bla (\cam{E})$ is the class of $\cam{E}$-modules
isomorphic to cokernels of maps in $\Sigma.$
\end{prop}

\begin{proof}
Note that the objects of $\Bla (\cam{E})$ are automatically
$\lev{E}$-modules, so that it makes sense to compare $\Bla
(\cam{E})$ and $\mathcal{E}(\lev{E},\Sigma)$.

Let us first show that $\Bla (\cam{E})=\mathcal{E}(\cam{E},\Sigma)$.
The proof follows arguments in \cite{FarberVogel92} and
\cite[Section 3]{RanickiSheiham06}; see also \cite[Section
6]{Ara04}. We will include most of the details for completeness.

First we show that the class $\mathcal{E}(\cam{E},\Sigma)$ is
exactly the class of Blanchfield $\cam{E}$-modules which are
finitely presented as left $\cam{E}$-modules. Since $\cam{E}$ is
hereditary, it suffices to show that, for a finitely presented
$\cam{E}$-module $M$, we have
$$\text{Tor}^{\cam{E}}_*(\cam{E}\Sigma^{-1}, M)=0 \iff
\text{Tor}^{\cam{E}}_*(\cb^d,M)=0.$$

Since $M$ is finitely presented there is an exact sequence
\begin{equation}
\label{exseforM}
\begin{CD}
0 @>>>P @>d>> Q @>>> M @>>> 0
\end{CD}
\end{equation}
with $P$ and $Q$ f.g. projective $\cam{E}$-modules. By \cite[Remark
3.4]{AraDicks07}, the map $1\otimes d\colon
\cam{E}\Sigma^{-1}\otimes_{\cam{E}} P\to
\cam{E}\Sigma^{-1}\otimes_{\cam{E}} Q$ is an isomorphism if and only
if the map $\epsilon (d):=1\otimes d\colon \cb^d
\otimes_{\cam{E}}P\to \cb^d\otimes_{\cam{E}} Q$ is an isomorphism.

For a module $X$, we use the canonical projective resolution of
$\cb^d$
\begin{equation*}
\begin{CD}
0 @>>> \bigoplus_{e\in E^1} p_{r(e)}\cam{E} @>{(e)}>> \cam{E} @>>>
\cb^d @>>> 0
\end{CD}
\end{equation*}
to compute the groups $\text{Tor}_*^{\cam{E}}(\cb^d,X)$. It follows
that $X$ is a Blanchfield $\cam{E}$-module if and only if the map
$\gamma_X\colon \bigoplus_{e\in E^1} p_{r(e)}X\to X$,  $\gamma
_X((p_{r(e)}x_e))=\sum ex_e$, is an isomorphism. Now the diagram in
the proof of \cite[Proposition 3.9(i)]{RanickiSheiham06} shows that
for the f.p. module $M$ with presentation (\ref{exseforM}), we have
that $\gamma _M$ is an isomorphism if and only if $\epsilon (d)$ is
an isomorphism. Hence, by the above comments, $M$ is a Blanchfield
module if and only $M$ is a $(\cam{E},\Sigma)$-torsion module.

To finish the proof that $\mathcal{E}(\cam{E},\Sigma)=\Bla
(\cam{E})$, we have to show that every f.g. Blanchfield
$\cam{E}$-module is finitely presented as $\cam{E}$-module. For this
part, we follow \cite[proof of Lemma 4.3]{FarberVogel92}.

Let $M$ be a f.g. Blanchfield $\cam{E}$-module.  Let $A$ be a
lattice in $M$ (Proposition \ref{typeLL}(2)), and consider the
$\cam{E}$-module endomorphism of the f.g. projective
$\cam{E}$-module $\cam{E}\otimes_{\cb^d} A$:
$$u\colon \cam{E}\otimes _{\cb^d} A\to \cam{E}\otimes _{\cb^d} A, \qquad u(\lambda \otimes a)=
\lambda \otimes a -\sum_{e\in E^1} \lambda e\otimes \ol{e} a ,$$
where $\lambda \in \cam{E}$ and $a\in A$. . Clearly $\epsilon
(u)=1$, and thus $\coker (u)\in \mathcal{E}(\cam{E},\Sigma)$. So the
previous argument gives that $\coker (u)\in \Bla (\cam{E})$. Let
$f\colon \cam{E}\otimes_{\cb^d} A\to M$ be the map given by
$f(\lambda \otimes a)=\lambda a$. Since $M$ is a Blanchfield module
we have $fu=0$, and thus there is a homomorphism $g\colon \coker
(u)\to M$ given by $g([\lambda \otimes a])=\lambda a$. The map
$\psi\colon A\to \coker (u)$, $\psi (a)=[1\otimes a]$ is
$\cam{\ol{E}}$-linear. Indeed we have, for $e'\in E^1$,
$$\ol{e'}\psi (a)=\ol{e'}[1\otimes a]=\ol{e'}[\sum _{e\in E^1}e\otimes \ol{e}a]=
\sum_{e\in E^1}\ol{e'}e[1\otimes \ol{e}a]=[1\otimes \ol{e'}a]=\psi
(\ol{e'}a).$$ We clearly have the identity $\iota =g \psi$, where
$\iota\colon A\to M$ denotes the inclusion. In particular $\psi$ is
injective and so $A$ is isomorphic with the lattice $\psi (A)$ of
$\coker (u)$. By Proposition \ref{typeLL}(2), the maps $1\otimes
\psi \colon \lev{E}\otimes _{\cam{\ol{E}}} A\to \coker (u)$ and
$1\otimes \iota\colon \lev{E}\otimes _{\cam{\ol{E}}}A\to M$ are both
isomorphisms, and clearly $1\otimes \iota =g(1\otimes \psi)$. It
follows that $g=(1\otimes \iota) (1\otimes \psi )^{-1} $ is an
isomorphism, so that in particular $M$ is finitely presented as a
$\cam{E}$-module. Moreover, this argument also shows the last
statement in the proposition.

Now we will show that $\mathcal{E}(\lev{E},\Sigma)=\Bla (\cam{E})$.
For $M\in \Bla (\cam{E})$ we have a projective resolution
\begin{equation*}
\begin{CD}
0 @>>> \cam{E}^n @>\sigma>> \cam{E}^n @>>> M @>>> 0 \, ,
\end{CD}
\end{equation*}
with $\sigma\in \Sigma$. Note that $ \sigma \colon \lev{E}^n \to
\lev{E}^n$ is also injective because the universal localization
$L(E)\to Q(E)=L(E)\Sigma^{-1}$ is injective. Thus we get a
resolution of $\lev{E}\otimes _{\cam{E}}M$:
\begin{equation}
\begin{CD}
0 @>>> L(E)^n @>\sigma>> L(E)^n @>>> L(E)\otimes _{\cam{E}}M @>>> 0.
\end{CD}
\end{equation}
Being $M$ an $\lev{E}$-module ,we get $L(E)\otimes _{\cam{E}}M\cong
M$, and thus $M\in \mathcal{E}(L(E),\Sigma)$.

Now it is straightforward to show that $\Bla
(\cam{E})=\mathcal{E}(\cam{E},\Sigma)$ satisfies (1)--(4) in
Definition \ref{def:catNee} for the pair $(L(E),\Sigma)$, hence we
get $\Bla (\cam{E})=\mathcal{E}(L(E),\Sigma)$, as desired.
\end{proof}

In our concluding result we compute the $K$-groups of the regular
algebra $Q(E)$. The Grothendieck group $K_0(Q(E))$ was computed in
\cite[Theorem 4.2]{AraBrustenga07}. We write ${\rm
Bla}_*(P(E))=K_*(\Bla (P(E)))$ for the $K$-groups of the exact
category $\Bla (P(E))$. As a preparation we compute $K_i(\rat{E})$.

\begin{lem}
\label{lem:Krat} Let $E$ be a finite quiver with $|E ^0|=d$. Then
there is a split exact sequence, for $i\ge 1$,
\begin{equation}
\label{eq:Krat}
\begin{CD}
0 @>>> K_i(\cam{E}) @>>> K_i(\cam{E}\Sigma^{-1}) @>>> {\rm
Bla}_{i-1}(\cam{E}) @>>> 0 \, ,
\end{CD}
\end{equation}
and so $K_i(\cam{E}\Sigma^{-1})=K_i(\rat{E})= K_i(\cb)^d\oplus {\rm
Bla}_{i-1}(\cam{E})$.
\end{lem}

\begin{proof}
Since $\cam{E}$ is hereditary, we can apply Theorem
\ref{theor:Nee-Ran} to the universal localization $\cam{E}\to
\cam{E}\Sigma^{-1}=\rat{E}$ to obtain an exact sequence in
nonnegative $K$-theory
$$\cdots \to K_{i}(\cam{E})\to K_{i}(\cam{E}\Sigma^{-1})\to
{\rm Bla}_{i-1}(\cam{E})\to K_{i-1}(\cam{E})\to \cdots .$$

We first show that the canonical embedding $\cb ^d\to \cam{E}$
induces an isomorphism $K_*(\cb^d)\to K_*(\cam{E})$ for $*\ge 0$.
This follows from \cite[Theorem 3.1]{Gersten74}, once we observe
that $P_\cb (E)[t]=P_{\cb [t]}(E)$ is regular coherent in the sense
of \cite{Gersten74}. The latter assertion follows from
\cite[Proposition 1.9 and Remark 1.10]{Gersten74}, by using
induction on the number of arrows of $E$, taking into account that
$P_A(E)\cong P_A(E')*_{A^d} P_A(E'')$, where $E'$ and $E''$ are
subquivers of $E$ with the same vertices and such that $E^1$ is the
disjoint union of $E'^1$ and $E''^1$. The basic case is the one in
which the quiver $E$ only has one arrow. If this arrow is a loop
then $P_{\cb [t]}(E)$ is clearly regular coherent because the
polynomial rings $\cb [t]$ and $\cb[t,s]$ are Noetherian regular
rings. If the arrow is not a loop then we get a triangular ring over
$\cb [t]$, and this is again Noetherian regular.

Now note that the isomorphism $K_i(\cam{E})\to K_i(\cb ^d)$, which
is induced by the augmentation map, factors through
$K_i(\cam{E}\Sigma^{-1})$, and so we see that the map
$K_i(\cam{E})\to K_i(\cam{E}\Sigma^{-1})$ has a retraction and, in
particular, it is injective. This shows the result.
\end{proof}

\begin{theor}
\label{theor:K1QE} Let $E$ be a finite quiver with $|E ^0|=d$. Then
$Q(E)$ is the universal localization of $\cam{\ol{E}}$ with respect
to the set of all monomorphisms between finitely generated
projective left $\cam{\ol{E}}$-modules whose cokernel is
finite-dimensional and does not contain nonzero projective modules.
Moreover we have, for $i\ge 1$,
$$
K_i(Q(E))  \cong K_i(L(E))\bigoplus {\rm Bla}_{i-1}(P(E)).$$
 In
particular
\begin{align*}
K_1(Q(E)) & \cong \coker (1-N_E\colon (\cb^{\times})^{(E_0\setminus \Si
(E))}\longrightarrow (\cb^{\times})^{(E_0)})\\
&  \bigoplus \ker (1-N_E\colon \Z^{(E_0\setminus \Si
(E))}\longrightarrow \Z^{(E_0)})\bigoplus {\rm Bla}_0(P(E))
\end{align*}
\end{theor}

\begin{proof}
Let $\Upsilon$ be the class of all monomorphisms between f.g.
projective $\cam{\ol{E}}$-modules whose cokernel is
finite-dimensional and does not contain nonzero projective modules.
Let $\Upsilon '$ be the class of monomorphisms between  f.g.
projective $\lev{E}$-modules induced by $\Upsilon$. Since the maps
$\nu _i$, for $i\in E^0\setminus \Si (E)$ (defined in the
Introduction), are in $\Upsilon$, we see that $\cam{\ol{E}}\Upsilon
^{-1}=\lev{E}\Upsilon '^{-1}$.

By Proposition \ref{typeLL}, we have that $\Bla (P(E))\cong
\fnp{\lev{E}}$ is exactly the class of cokernels of maps in
$\Upsilon '$.  Since $\Bla (\cam{E})=\mathcal{E}(\lev{E},\Sigma )$
by Proposition \ref{torsions}, it follows that
$$Q(E)=\lev{E}\Sigma^{-1}=\lev{E}\Upsilon '^{-1}=\cam{\ol{E}}\Upsilon
^{-1}.$$
This shows the first part of the theorem.

Since both $\cam{E}$  and $\lev{E}$ are hereditary, we can apply
Theorem \ref{theor:Nee-Ran} to the two universal localizations
$\cam{E}\to \cam{E}\Sigma^{-1}$ and $\lev{E}\to \lev{E}\Sigma
^{-1}=Q(E)$. Comparison of both localization sequences gives, taking
into account Lemma \ref{lem:Krat}, the following commutative diagram
of exact sequences, for $i\ge 1$:
\begin{equation}
\begin{CD}
\label{eq:qqq} 0 @>>> K_i(P(E)) @>>> K_i(\rat{E}) @>>> {\rm Bla}_{i-1}(P(E)) @>>> 0\\
 & & @VVV  @VVV  @VV=V  \\
 & & K_i(L(E)) @>>> K_i(Q(E)) @>>> {\rm Bla}_{i-1}(P(E))
\end{CD}
\end{equation}
It follows that the map $K_i(Q(E)) \to {\rm Bla}_{i-1}(P(E))$  is
surjective, and so we get a short exact sequence, for $i\ge 1$,
\begin{equation}
\label{eq:KQE}
\begin{CD} 0 @>>> K_i(L(E))  @>>> K_i(Q(E))@>>> {\rm
Bla}_{i-1}(P(E)) @>>> 0
\end{CD}
\end{equation}
Since the exact sequence (\ref{eq:Krat}) splits, so does the exact
sequence (\ref{eq:KQE}), by (\ref{eq:qqq}). The formula for
$K_1(Q(E))$ follows now from \cite{AraBrustengaCortinas09}.
\end{proof}

\bibliographystyle{plain}
\bibliography{biblioTeX}
\end{document}